\documentclass[reprint,aps,pre,twocolumn]{revtex4-2} % two-column body, wide abstract
\usepackage{graphicx}
\usepackage{amsmath}
\usepackage{amssymb}
\usepackage{booktabs}
\usepackage{hyperref}

\newcommand\diff{\mathrm{d}}

\newcommand{\mathsfi}[1]{\mathsf{#1}}
\renewcommand{\vec}[1]{\mathbf{#1}}
\newcommand{\vect}[1]{\boldsymbol{#1}}
\newcommand{\T}{\mathrm{T}}
\renewcommand{\tensor}[1]{\boldsymbol{\mathsfi{#1}}}
\renewcommand{\L}{\mathcal{L}}

\usepackage{bm}

\begin{document}

\title{Meshless solutions of PDE inverse problems on irregular geometries}
\author{James V. Roggeveen}
\email{roggeveen@seas.harvard.edu}
\affiliation{School of Engineering and Applied Sciences, Harvard University, Cambridge MA 02138}
\author{Michael P. Brenner}
\email{brenner@seas.harvard.edu}
\affiliation{School of Engineering and Applied Sciences, Harvard University, Cambridge MA 02138}
\affiliation{Department of Physics, Harvard University, Cambridge MA 02138}

\begin{abstract}
Solving inverse and optimization problems over solutions of nonlinear partial differential equations (PDEs) on complex spatial domains is a long-standing challenge. 
Here we introduce a method that parameterizes the solution using spectral bases on arbitrary spatiotemporal domains, whereby the basis is defined on a hyperrectangle containing the true domain. We find the coefficients of the basis expansion by solving an optimization problem whereby both the equations, the boundary conditions and any optimization targets are enforced by a loss function, building on a key
idea from Physics-Informed Neural Networks (PINNs). Since the representation of the function natively has exponential convergence, so does the solution of the optimization problem, as long as it can be solved efficiently. We find empirically that the optimization protocols developed for machine learning find solutions with exponential convergence on a wide range of equations. The method naturally allows for the incorporation of data assimilation by including additional terms in the loss function, and for the efficient solution of optimization problems over the PDE solutions.
\end{abstract}

\maketitle

Solving inverse problems on the solutions of partial differential equations (PDEs) has emerged as a core challenge in many fields of science, mathematics, and engineering. Applications range from data assimilation\cite{law2015data}, whereby a numerical solution to a PDE is fit to experimental data, to control problems,
where some portion of the solution, such as the boundary shape, is optimized to achieve a specified target objective \cite{Zhou2024, Alhashim2025}.

Such inverse problems are often ill-posed, with insufficient boundary information for a traditional solver, noisy data, and hard geometric constraints. Typically, optimization targets are met by finding gradients of the PDE solution with respect to the parameters. While such problems have been traditionally addressed using adjoint methods \cite{giles2002adjoint}, the recent proliferation of automatic differentiation tools has enabled the construction of differentiable PDE solvers \citep{Alhashim2025}, whose gradients can be used to replace adjoint solves in an optimization loop. However, both adjoint and differentiable methods require many iterations of a PDE solver to converge on the optimal solution. Further, deriving adjoint equations in arbitrarily complex domains can be immensely challenging and applying differentiable solvers requires reimplementing traditional finite-element \cite{LeVeque2007} or spectral \citep{Trefethen2000,Shen2011} methods within a differentiable code.

Several alternative approaches have recently emerged relying on neural networks to represent a solution to a family of PDEs, including Fourier neural operators (FNOs) \cite{Li2021} and Physics-Informed Neural Networks (PINNs) \citep{Raissi2019,WangS2023}. PINNs find solutions to PDEs by representing the solution with a deep neural network, with weights chosen to minimize the PDE residuals across a set of collocation points spanning the domain. The solution to this optimization problem leverages automatic differentiation for computing the derivatives in the PDE residuals and for minimizing the residual loss. 

PINNs have been demonstrated to be a remarkably flexible method for solving PDEs \citep{Raissi2019,WangS2023}, particularly for inverse problems and data assimilation \citep{Wang2025}. By representing the solution as a neural network, PINNs are naturally grid-free and do not require meshing, lending themselves to solving problems on complex geometries. Additionally, it is simple to augment the global loss function over residuals with additional terms that enforce data-assimilation matching constraints, such as enforcing that the solution obeys a set of measurements. In these cases, the total loss function $\L$ is written as the sum
\begin{equation}
  \L = \L_{\mathrm{pde}} + \L_{\mathrm{boundary}} + \L_{\mathrm{data}},
  \label{eq:PINNLoss}
\end{equation}
where $\L_{\mathrm{pde}}$ represents the loss due to the PDE residuals at a set of collocation points, $\L_{\mathrm{boundary}}$ the residuals of the boundary conditions at a set of boundary collocation points, and $\L_{\mathrm{data}}$ the squared error between observed data and the solution. 

Yet, the flexibility of PINNs is not without costs. While neural networks are theoretically universal function approximators \citep{Cybenko1989,Hornik1989}, in practice PINNs are constrained by spectral biases \citep{Cao2019,Rahaman2019} that result from truncation error of a finite neural network architecture. Furthermore, PINNs lack the convergence guarantees of classical numerical analysis. To compensate, they are typically extremely overparameterized, with networks ranging from 60,000 to over 100,000 individual parameters, leading to computationally expensive gradient calculations and unstable training dynamics \citep{Wang2021}. PINNs also can struggle with higher-order PDEs, where the computational cost of computing derivatives scales exponentially with order \citep{Karnakov2023}. Finally, PINNs offer little in the way of immediate physical interpretability, as the weights of a given network do not directly map onto standard functions.

Another recent approach has sought to combine the flexibility of the PINNs global loss formulation (Eqn.~(\ref{eq:PINNLoss})) with standard PDE techniques by representing functions as fixed nodes on a grid and using finite-difference stencils to compute gradients for the purpose of PDE residuals \cite{Karnakov2023}. This results in dramatic computational speedups, while maintaining many of the PINNs' advantages in terms of data assimilation and flexible loss functions. 
However, this approach still requires grid or mesh generation, complicating the application to very complex domains.

\begin{figure}
  \centering
  \includegraphics[width=\linewidth]{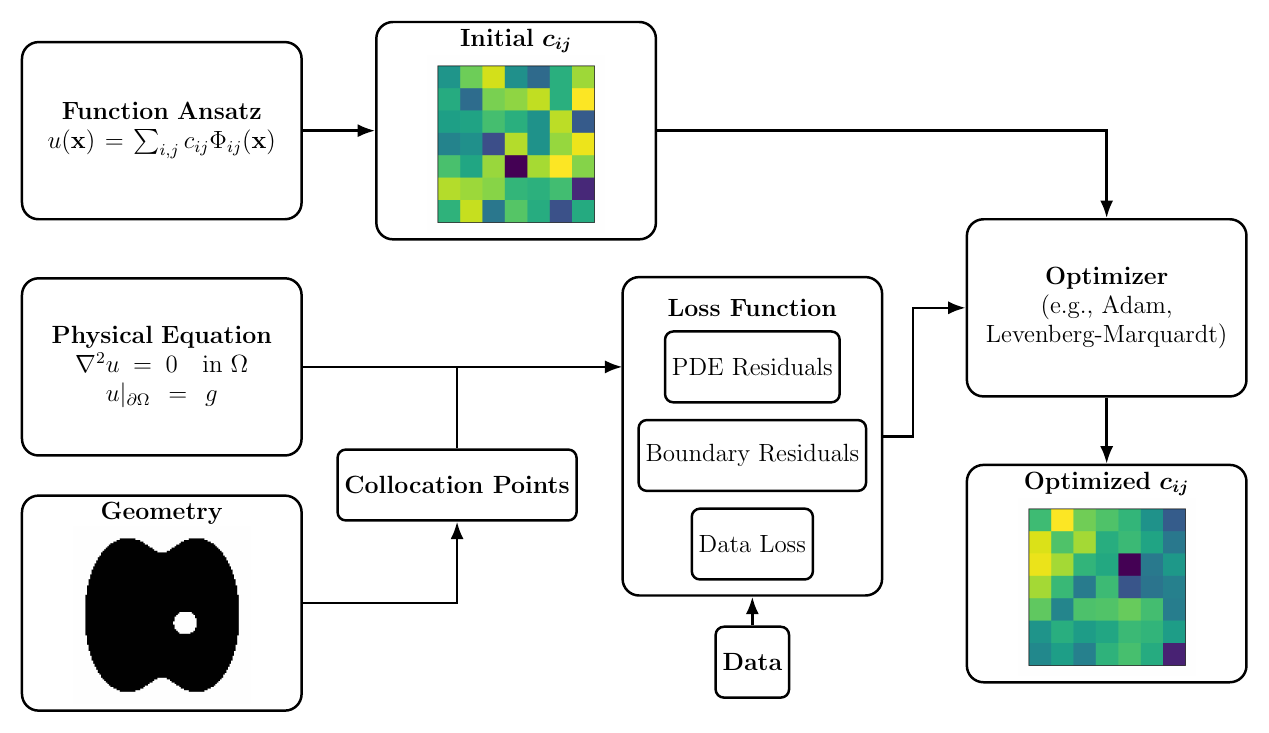}
  \caption{Schematic view of basis fitting scheme. We assume that any function can be written as a tensor product of basis terms with some set of coefficients, which are not initially known. We define a loss function over some geometry by computing the residuals of some physical equation and boundary conditions at a series of collocation points. The loss function can also include terms that require the solution to match data. Find a set of optimized coefficients $c_{ij}$ that minimize this loss function using a differentiable code and gradient-based nonlinear optimizer.}
  \label{fig:schematic}
\end{figure}

%Achieving better convergence in error requires representing the solution using spectral bases. One such approach is to use radial basis functions defined at a series of collocation points to smoothly represent solutions on arbitrary geometries \citep{Fornberg2015,Flyer2011}. One major downside to RBF methods is that computing the derivatives of the functions can become expensive as the number of nodes increases. Further, traditional RBF methods become ill-conditioned as the functions are smoothed out to provide better approximations and generally rely on a time-stepping a approach to solve unsteady problems. 

Here we build on these ideas and demonstrate a new method for solving PDEs and associated inverse problems in complex domains (Fig.~\ref{fig:schematic}). We use the global residual loss (Eqn.~(\ref{eq:PINNLoss})) but represent the solution using a spectral basis over a hyperrectangular domain that {\sl contains} the complex boundary shape, with the boundary condition itself enforced via the loss function.
We demonstrate through a series of examples that using the nonlinear optimizers developed by the machine learning community, we can solve a series of linear and nonlinear forward and inverse problems in multiple dimensions in complex domains, with empirically observed exponential convergence and an order of magnitude fewer parameters than neural network-based methods.

\section{Problem Formulation}

We represent the solution $u$ over a hyper-rectangular domain with
coordinates $\vec{x} = (x_1,...,x_d)^\T$ with the expansion
\begin{equation}
  u = \tensor{C}:\vect{\Phi}(\vec{x}),
  \label{eq:sol_func}
\end{equation}
where $\vect{\Phi}(\vec{x}) = \vect{\phi}_1(x_1)\dots\vect{\phi}_d(x_d)$ is the tensor product of spectral basis vectors $\vect{\phi}_i$, $\tensor{C}$ is the coefficient tensor of rank $d$ and the double dot indicates a Frobenius inner product. Note that both the number and type of basis functions in each $\vect{\phi}_i$ need not be the same. Thus, we can mix and match functions, such as Fourier series, Chebyshev polynomials, or Legendre polynomials, and degrees of accuracy as needed for a given problem. It is important to note in our definition of a solution that we include time $t$ as a coordinate whose solution we encode with a spectral basis. Thus, a transient problem defined over a three-dimensional volume would be represented by the coordinate vector $\vec{x} = (t,x,y,z)^\T$. As such, we handle time in the same way as any spatial coordinates rather than with an explicit time-stepping algorithm.

We find the optimal set of coefficients $\tensor{C}$ that minimizes Eqn.~(\ref{eq:PINNLoss}). While this is in general nonconvex, recent advances in gradient-based optimization methods (e.g. \cite{zhang2018improved}) make this a tractable optimization problem to solve, albeit with the loss of uniqueness guarantees on the solution. 

This approach has a number of benefits. A given field $u$ and its derivatives are defined everywhere in the domain without relying on fixed grids or interpolations. Since these derivatives are known analytically, the derivative of any function $u$ can be written with respect to derivatives of the tensor $\vect{\Phi}$, computed in linear time for any order without relying on automatic differentiation. For a fixed set of collocation points, the necessary derivative tensors $\partial^j/\partial x_i^j\vect{\Phi}$ may be precomputed before beginning the optimization, rendering any residual calculation to be a fast matrix multiplication. 
When the domain is non-rectangular, we use the Fourier extension theory to embed the computational domain inside of a hyper-rectangle upon which our basis functions are defined. Such an embedding has been shown \citep{Huybrechs2010} to still maintain favorable convergence properties at the cost of strict orthogonality. While the basis functions in these cases are not strictly orthogonal, they are an order of magnitude more parameter efficient than a neural network, leading to faster convergence in optimization. 

\section{Results}
\begin{figure}
  \centering
  \includegraphics[width=\linewidth]{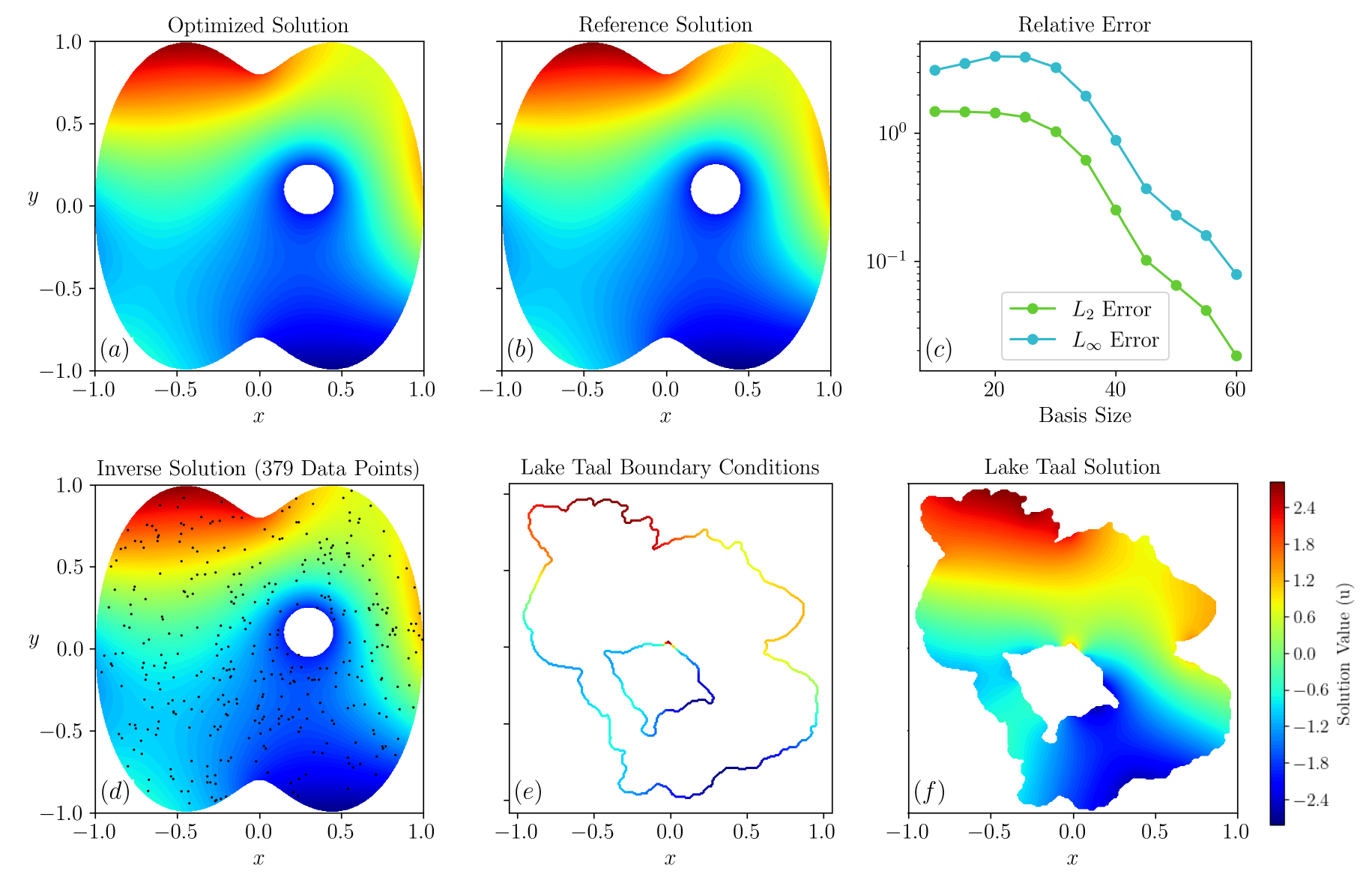}
  \caption{(a) Solution to the Laplace equation in an irregular domain with a hole as found by fitting polynomial coefficients to the equation residuals with $N=60$ coefficients in both $x$ and $y$. (b) The same system solved using COMSOL with an extremely high mesh density. (c) Convergence of the basis fitting approach to the COMSOL solution as the number of basis modes is increased. (d) Solution of Laplace's equation when the boundary condition is not provided but 379 data samples from the interior of the domain (shown as the black dots) included in the loss with $N=45$ coefficients. (e) Boundary condition for the Laplace equation on a very irregular domain described by the perimeter of Lake Taal in the Philippines. (f) Solution of the Laplace equation on the surface of Lake Taal with $N=100$ basis modes in both $x$ and $y$.}
  \label{fig:laplace}
\end{figure}
To test our basis fitting scheme at solving PDEs, we solve a series of problems with $d=2$ dimensions where a high-accuracy numerical solutions, $u_{ref}$, is known. We quantify error by evaluating our approximate solution $u$ at the fixed regular grid points $\vec{x}_{ij} = (x_{1,i},x_{2,j})$ on which the numerical solution is defined. We use both a discrete $L_2$ and $L_\infty$ errors.

\subsection{Laplace Equation}
We first solve the Laplace equation $\nabla^2 u =0$ in two dimensions.
To illustrate the ability to compute on complex domains, we choose a domain
 $\Omega$ as a peanut shape embedded in the rectangle $[-1,1]\times[-1,1]$, with a hole in it (Fig.~\ref{fig:laplace}(a)).
 The outer boundary $\partial\Omega_1$ is defined by a polar function $r(\theta)$ with a circular cut out, with boundary conditions $u(\partial \Omega)$, whose precise form is given in Methods.

We solve this equation by fitting $\tensor{C}$ to tensor product of Chebyshev polynomials defined in both $x$ and $y$ using a second-order Gauss-Newton optimizer running the Dogleg algorithm.
Fig.~\ref{fig:laplace}(a) shows the result for $N=60$ in each dimension, totaling 3721 coefficients (the coefficients are zero-indexed with respect to $N$). Fig.~\ref{fig:laplace}(b) shows the same problem solved using a finite element method implemented in COMSOL 6.3 with the highest automatic mesh density, which we take as our reference solution. Both the $L_2$ and $L_\infty$ error (Fig.~\ref{fig:laplace}(c)) converges exponentially with increasing basis modes $N$.

While our method can arrive at this forward solution to Laplace's equation in this simple peanut geometry, it does so with several orders of magnitude less accuracy and speed than more specialized methods, such as rational approximation \cite{Nakatsukasa2018}. For reference, the COMSOL numerical solution was obtained in less than a second compared to over six minutes for the solution shown in Fig.~\ref{fig:laplace}(a) when run on a laptop's CPU. However, what this method lacks in speed is made up for in flexibility. As an example, we again solve Laplace's problem on the same geometry in Fig.~\ref{fig:laplace}(d), but instead of using any information about the boundary condition we only enforce the equation residual and a matching condition to 379 randomly sampled points from the reference solution. This replacement requires no code changes to the solver. We find that in this case with $N=45$ modes in both directions a representation of the solution with an $L_2$ error of 0.02 and an $L_\infty$ error of 0.08, which is comparable to the solution with enforced boundary conditions. We demonstrate how the number of sample points affects the accuracy of the solver in Fig.~\ref{fig:laplaceinverseerror}.

To highlight the flexibility of our framework to solving problems on extremely complex geometries, we also solve the Laplace equation on a boundary defined by the shape of Lake Taal, embedded in a unit square. We use the same conditions for the outer boundary as the previous problem (Fig.~\ref{fig:laplace}(e)) applied now to both the outer and inner boundary on the lake. 
Switching to this new domain involved changing only a few lines of code to define the boundary from a PNG mask rather than an analytic equation.
Fig.~\ref{fig:laplace}(f) shows the solution for $N=100$ basis modes. The basis fitting method is extremely adaptable for solving problems in complex geometries for which finding a suitable mesh would be challenging.

\subsection{Time-dependent stiff equations}
\begin{figure}
  \centering
  \includegraphics[width=\linewidth]{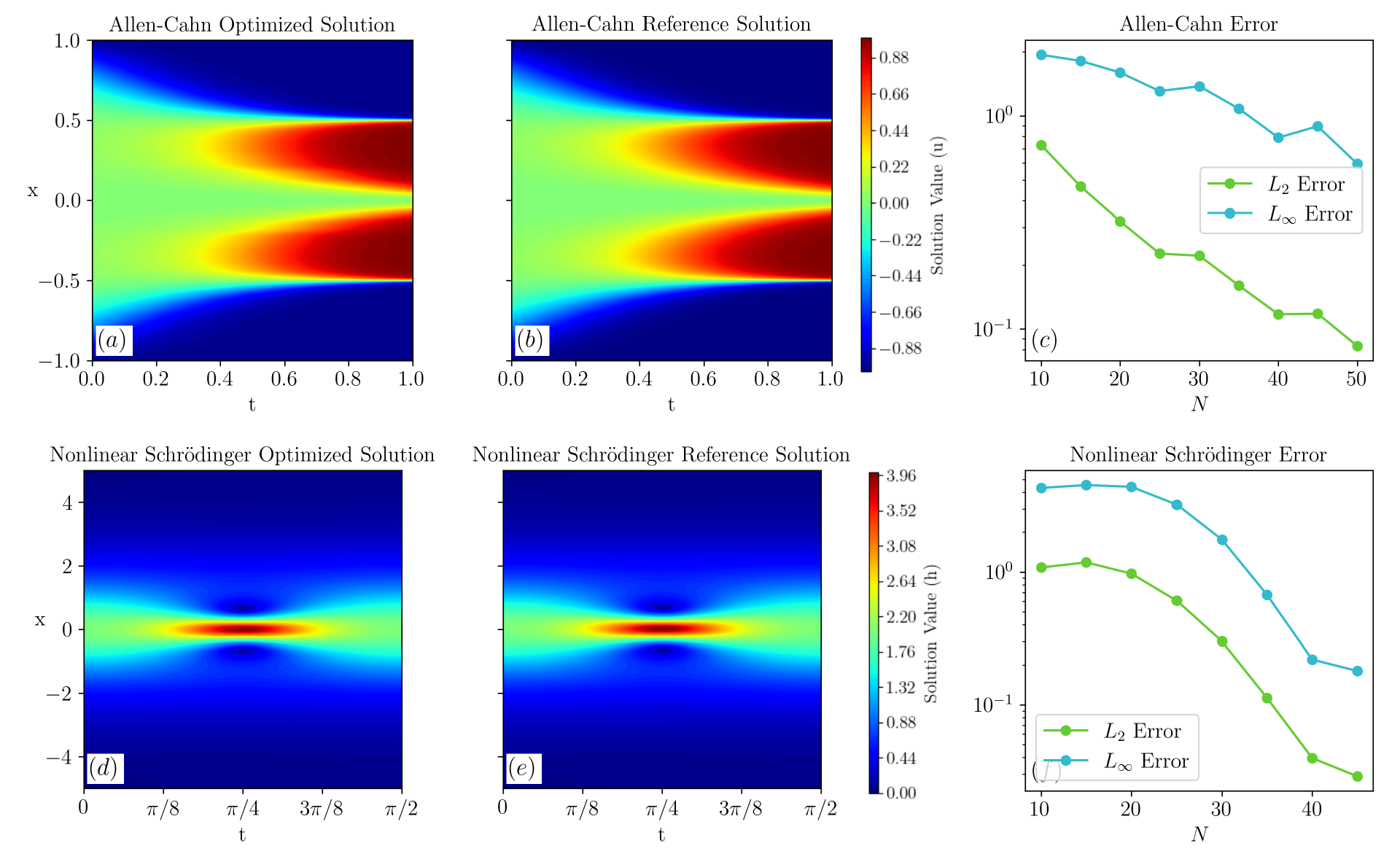}
  \caption{Solutions to stiff PDEs shown versus reference solutions. (a) shows a polynomial fit to the Allen-Cahn equations, with 200 basis modes in $x$ and 10 modes in $t$. (b) shows a reference solution from \cite{WangS2023}, while (c) shows the convergence of the solution to the reference as the number of basis modes in both $x$ and $t$ are increased. (d) similarly shows a polynomial fit to the nonlinear Schr\"odinger equations with $N=45$ in both $x$ and $t$ compared to a reference solution (e) from \cite{Raissi2019}. (f) shows the convergence of the polynomial fit to the reference as a function of the number of basis modes.}
  \label{fig:stiff_eqs}
\end{figure}
We now turn to time dependent problems where we demonstrate that not only does our method reliably find the correct solution, but empirical exponential convergence extends to the temporal domain. As illustrations, we explore two problems in one spatial dimension $x$ and time $t$ - the resulting domains are rectangles in $x-t$ space and do not include any non-regular boundaries. 

\subsubsection{Allen-Cahn equations} We first examine solutions to the Allen-Cahn equations, given by 
\begin{gather}
  \frac{\partial u}{\partial t} - 0.0001\frac{\partial^2 u}{\partial x^2} + 5u^3 - 5u = 0,\, t \in [0, 1], \, x \in [-1, 1], \label{eq:AC}
\end{gather}
with boundary condition $u(t, -1) = u(t,1)$, $\frac{\partial u}{\partial x}(t, -1) = \frac{\partial u}{\partial x}(t, 1)$ and initial conditions $u(x, 0) = x^2 \cos(\pi x)$.This equation is very stiff due to the small diffusion constant multiplying the spatial derivatives. We take the reference solution as that provided by \cite{WangS2023}. 
We choose to use a Fourier basis defined on $[-1,1]$ in the spatial domain, naturally satisfies the periodic boundary conditions on $u$ and its derivative. Since both Eqn.~(\ref{eq:AC}) and the initial condition are even in $x$, we restrict the Fourier basis to cosine modes. In time, we map the domain $[0,1]$ to $[-1,1]$ using an affine transformation and represent the function using Legendre polynomials. 
Fig.~\ref{fig:stiff_eqs}(a) shows the result of our method for $N_x=200$ and $N_t = 10$, with excellent agreement with the reference solution (Fig.~\ref{fig:stiff_eqs}(b)).  Fig.~\ref{fig:stiff_eqs}(c) shows evidence of exponential convergence even when the problem involves time. The inherent stiffness of the solution requires using a much higher number of basis modes in $x$ than in $t$ to achieve a qualitatively good agreement with the reference.

\subsubsection{Nonlinear Schr\"odinger Equations}
We next solve the nonlinear Schr\"odinger equations, given in terms of a complex-valued function $g$ as
\begin{gather*}
  i \frac{\partial g}{\partial t} + 0.5 \frac{\partial^2 g}{\partial x^2} + |g|^2 g = 0, \quad x \in [-5, 5], \quad t \in [0, \pi/2], 
\end{gather*}
with initial condition $g(0, x) = 2 \operatorname{sech}(x)$, and boundary conditions $g(t, -5) = g(t, 5)$ and $\frac{\partial g}{\partial x}(t, -5) = \frac{\partial g}{\partial x}(t, 5)$. 
We define the real and imaginary components as $g=u+iv$ and solve the resulting coupled nonlinear PDE. We again choose Fourier basis restricted to cosine modes since the solution is even and use a Chebyshev polynomial basis in $t$. The physical domain $[-5,5]\times[0,\pi/2]$ is mapped to $[-1,1]\times[-1,1]$ by an affine transformation in both variables. 

Fig.~\ref{fig:stiff_eqs}(d) shows the complex modulus $h = \sqrt{u^2 + v^2}$ of our optimized solution for $N=45$, compared to a reference numerical solution provided by \cite{Raissi2019} (Fig.~\ref{fig:stiff_eqs}(e)). We again analyze the convergence of our solution to the reference as a function of $N$ (Fig.~\ref{fig:stiff_eqs}(f)), with apparent exponential convergence. Here there is some evidence of saturation at the highest mode. This saturation may be a result of difficulties in optimizing a large parameter space over a fixed number of collocation points, as discussed more in Fig.~\ref{fig:laplace_analysis}. 

\subsection{Wave equation in a peanut}
\begin{figure}
  \centering
  \includegraphics[width=\linewidth]{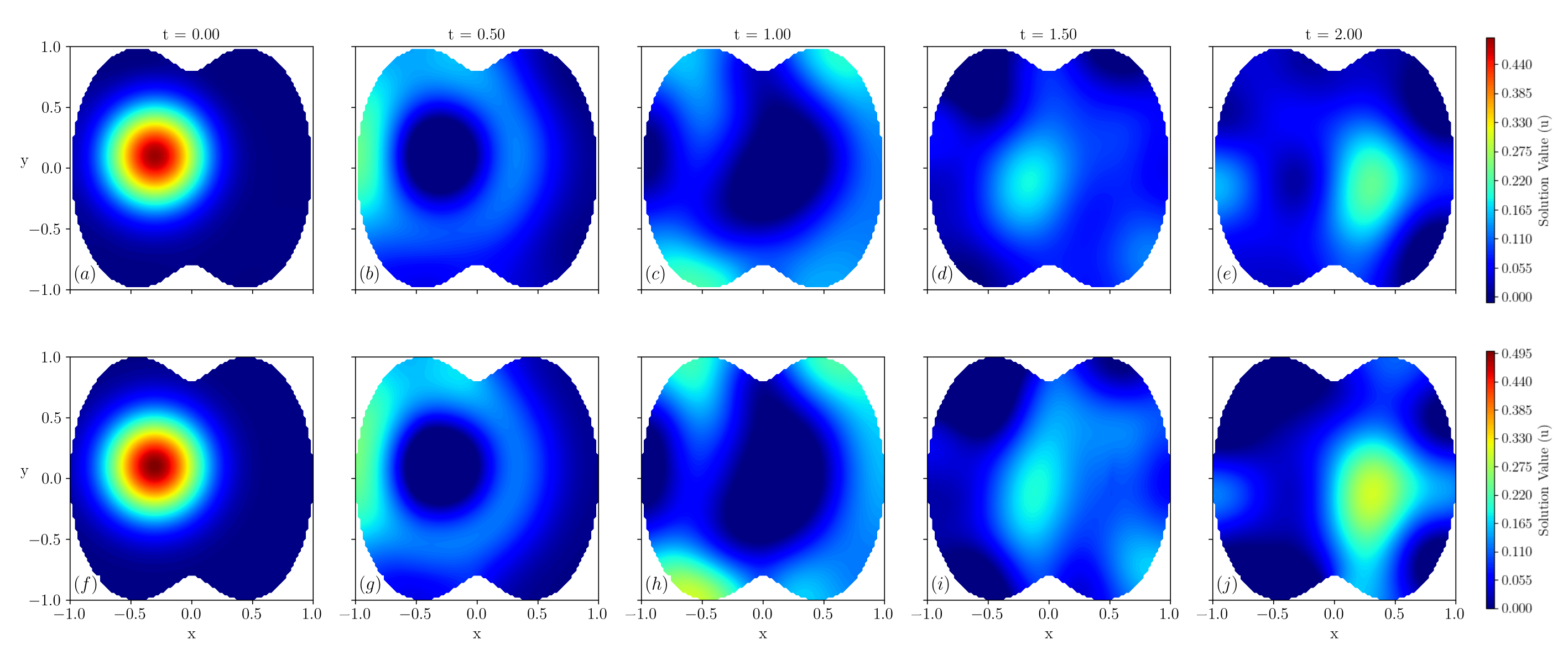}
  \caption{Solution to the wave equation on an irregular domain at several time steps evaluated using a polynomial fit (top row) with $N=12$ polynomials in each dimension of $x$, $y$, and $t$, and COMSOL (bottom row) solved using a finite element method with extremely fine mesh density.}
  \label{fig:wave}
\end{figure}
We return to the peanut domain $\Omega$ defined in the Laplace problem (ignoring the circular intrusion), and seek solutions to the time dependent wave equation: 
\begin{gather*}
	\frac{\partial^2 u}{\partial t^2} - \left( \frac{\partial^2 u}{\partial x^2} + \frac{\partial^2 u}{\partial y^2} \right) = 0, \, \forall (x,y) \in \Omega, \, t \in [0,2],
\end{gather*}
with Neumann boundary conditions $\nabla u \cdot \vec{n} = 0$ on $\partial\Omega$ 
and a Gaussian pulse initial condition $u(x, y, 0) = 0.5 \exp\left(-\frac{(x+0.3)^2 + (y-0.1)^2}{0.16}\right)$.
We solve the problem with $N=12$ modes of Legendre polynomials in each of $x$, $y$, and $t$. 

The top row of Fig.~\ref{fig:wave} shows snapshots of our solutions, while the bottom row shows the output of a COMSOL using the same domain and initial conditions at the maximum automatic mesh resolution. Our method, even with a small number of modes, finds topologically similar solutions in space and time, with an $L_2$ error of $0.034$. 
%The magnitude of our solution is generally slightly lower than that of the reference, which can be understood by examining the lower magnitude of the initial pulse in the first column. 
Using only $N=12$ Legendre modes cannot completely capture the exact Gaussian initial condition, leading to less overall energy in the system and smaller amplitudes in the final solution. This discrepancy could be ameliorated by increasing the size of the spatial basis modes.

\subsection{Diffusion on the surface of a sphere}
\begin{figure}
  \centering
  \includegraphics[width=\linewidth]{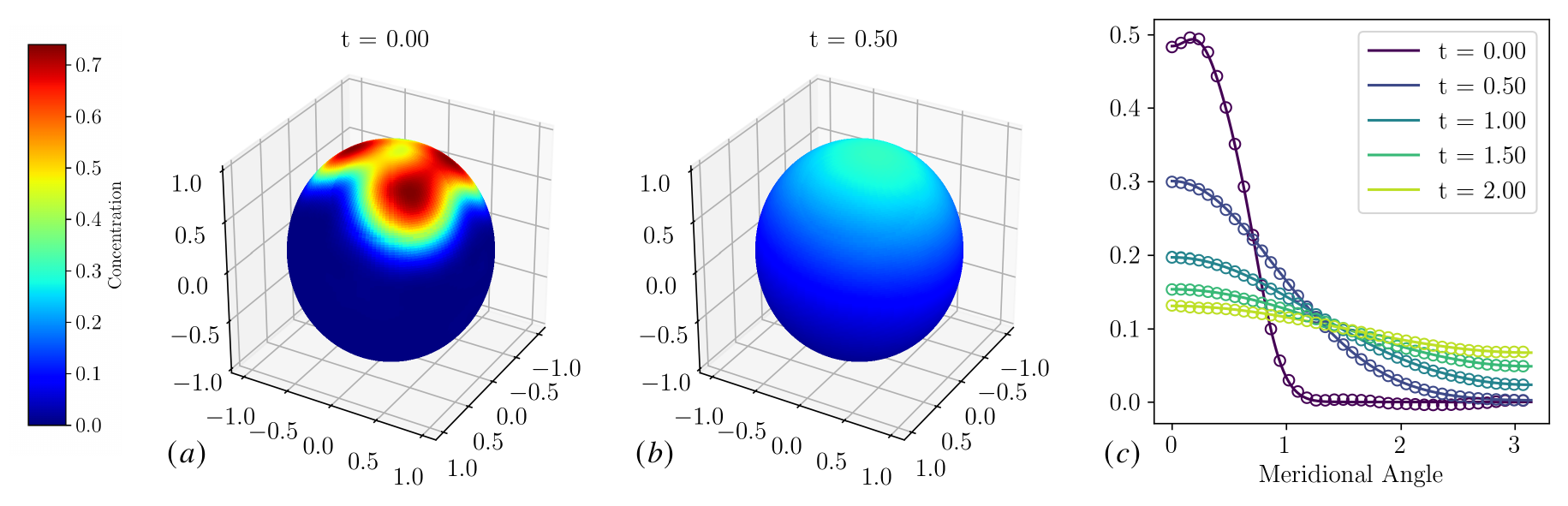}
  \caption{Diffusion equation solved on the surface of a sphere. The surface was embedded in 3D space and the surface diffusion equation was solved with $N=10$ Chebyshev modes in each of $x$, $y$, $z$, and $t$. Simulations were run to $t=2.0$ with $D = 0.5$. (a) and (b) show the concentration field on the sphere from the non-axisymmetric initial condition at $t=0$ to the increasingly axisymmetric profile at $t=0.5$. (c) compares the our simulation results in the open circles with the exact analytic result (solid line) expressed in spherical harmonic functions with 20 modes. Results are shown at several time points in the simulation. }
  \label{fig:spherediff}
\end{figure}
We next demonstrate solutions of problems on more complex manifolds by embedding them within a higher dimensional space. In particular, we seek a solution to the diffusion equation confined to the surface of a unit sphere, with governing equations
\begin{gather}
  \frac{\partial c}{\partial t} = D\nabla_S^2 c=D \left[ \nabla^2 c - \frac{\partial^2 c}{\partial n^2} - 2 \frac{\partial c}{\partial n} \right]
\end{gather}
and initial conditions $c(\vec{x},0) = c_0(\vec{x}) = \sum_{i=1}^3S(\vec{x},\vec{x}_{c,i})$, where 
$S(\vec{x},\vec{x}_c) = 0.5(\tanh[10(\vec{x}\boldsymbol{\cdot}\vec{x}_c - 0.95)] + 1)$ and
$\vec{x}$ are points on the surface of the unit sphere. 
We take the initial condition to be centered at the points 
$\vec{x}_{c,1} = (-0.40,0.40,0.82)^\T$, 
$\vec{x}_{c,2}=(0.40,0.40,0.82)^\T$ and $\vec{x}_{c,3} = (0.00,-0.40,0.92)^\T$.

We solve the diffusion problem using with $D=0.5$ for $t\in[0,2]$, choosing $N=10$ Chebyshev modes in each of the four dimensions. Fig.~\ref{fig:spherediff} (a) and (b) show snapshots of the concentration field at $t=0$ and $t=0.5$.

We can examine the accuracy of this solution compared against an analytic solution in terms of spherical harmonics. Fig.~\ref{fig:spherediff}(c) compares our solution with a sum of 20 spherical harmonic modes at different time points along a meridian, with excellent agreement over all times. 
We note that while we chose this particular example to have a known analytic solution to investigate solution accuracy, the  method can solve more complex equations on more complex surfaces. This has proven difficult for finite element methods \cite{Landsberg2010} and impossible for typical spectral methods. 

\subsection{Inversion of ice sheet viscosity}
\begin{figure}
  \centering
  \includegraphics[width=\linewidth]{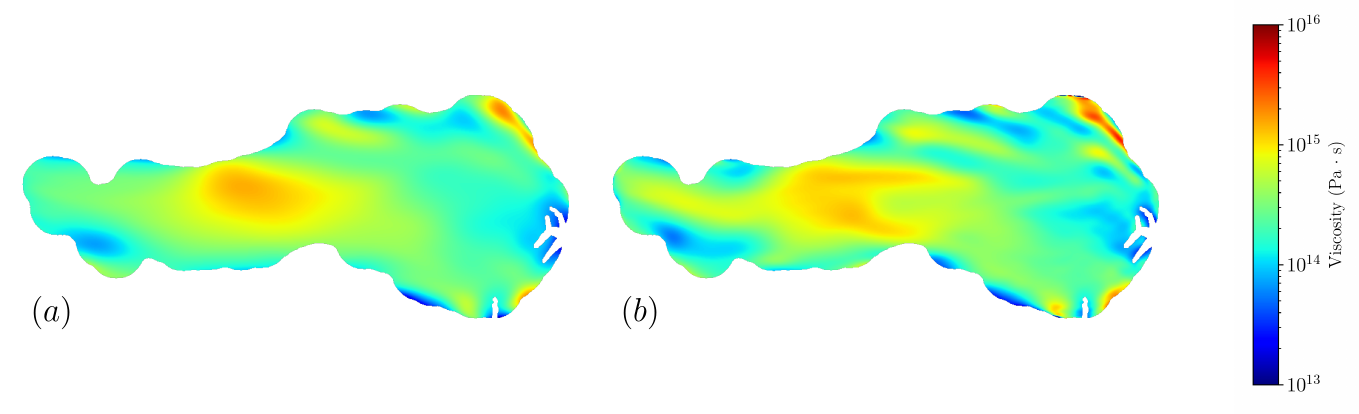}
  \caption{Viscosity inversion through the shallow shelf approximation on the Amery ice sheet using data published in \cite{Wang2025}. (a) uses $N=25$ Chebyshev modes in each coordinate to express the field while (b) uses $N=50$ modes. To fit the viscosity, the coefficients were used to find the log of the viscosity to help account for the differences in scale across the ice sheet.}
  \label{fig:ameryinversion}
\end{figure}
We now turn to a real-world problem of finding the viscosity profiles of Antarctic ice sheets given observational data, recently investigated \cite{Wang2025} with PINNs. These authors solve the steady shallow shelf approximation equations in two dimensions,
\begin{multline}
\frac{\partial}{\partial x}\left(4\mu h \frac{\partial u}{\partial x} + 2\mu h \frac{\partial v}{\partial y}\right) + \frac{\partial}{\partial y}\left(\mu h \frac{\partial u}{\partial y} + \mu h \frac{\partial v}{\partial x}\right) \\= \rho_i g \left(1 - \frac{\rho_i}{\rho_w}\right) h \frac{\partial h}{\partial x},
\label{eq:ssa_eq1}
\end{multline}
\begin{multline}
\frac{\partial}{\partial y}\left(4\mu h \frac{\partial v}{\partial y} + 2\mu h \frac{\partial u}{\partial x}\right) + \frac{\partial}{\partial x}\left(\mu h \frac{\partial u}{\partial y} + \mu h \frac{\partial v}{\partial x}\right) \\= \rho_i g \left(1 - \frac{\rho_i}{\rho_w}\right) h \frac{\partial h}{\partial y}.
\label{eq:ssa_eq_2}
\end{multline}
Here, $h$ is the thickness of the ice sheet, $(u,v)$ is the velocity of the sheet in $(x,y)$, $\mu$ is a spatially varying viscosity, and $\rho_i (\rho_w)$ are the density of ice (water), taken to be 917 and 1030 $\mathrm{kg}/\mathrm{m}^3$. In addition to these equations, there is a force balance condition applied at the calving front, 
\begin{multline}
2\mu h \left(2\frac{\partial u}{\partial x} + \frac{\partial v}{\partial y}\right) n_x + \mu h \left(\frac{\partial u}{\partial y} + \frac{\partial v}{\partial x}\right) n_y \\ = \frac{1}{2}\rho_i g \left(1 - \frac{\rho_i}{\rho_w}\right) h^2 n_x,
\label{eq:ssa_bc_1}
\end{multline}
\begin{multline}
\mu h \left(\frac{\partial u}{\partial y} + \frac{\partial v}{\partial x}\right) n_x + 2\mu h \left(\frac{\partial u}{\partial x} + 2\frac{\partial v}{\partial y}\right) n_y \\= \frac{1}{2}\rho_i g \left(1 - \frac{\rho_i}{\rho_w}\right) h^2 n_y.
\label{eq:ssa_bc_2}
\end{multline}
The calving front has a normal vector $\vec{n}=(n_x,n_y)^T$. 

In this particular problem, observational satellite data is available for the sheet thickness $h$ along with the sheet velocity $u$ and $v$. One difficulty is that the grids on which these data are defined are not the same and the ice sheet itself has a very irregular shape, making grid-based solvers impractical.

Our goal is to find representations of the fields $h$, $u$, $v$, and $\mu$ that satisfy the equations while agreeing with the data. We thus include a data loss component $\L_{\mathrm{data}}$ equal to the sum of the squared residuals between the fitted functional forms and a sample of the experimental data. We non-dimensionalize the equations to equalize the scales in the problem and use an affine transformation to embed the ice sheet in $[-1,1]\times[-1,1]$. Finally, given that the viscosity is expected to vary over several orders of magnitude, we solve for $\log \mu = \tensor{C}:\vect{\Phi}$. 

We solve the problem using a Chebyshev basis in both directions and show the resultant dimensional viscosity fields for $N=25$ modes in Fig.~\ref{fig:ameryinversion}(a), and with $N=50$ modes in Fig.~\ref{fig:ameryinversion}(b). The topology is similar in structure and magnitude to that predicted with PINNs \cite{Wang2025}, with an increase in structures and detail with increasing numbers of basis modes. Verification of this method on a known COMSOL solution from \cite{Wang2025} is demonstrated in Fig.~\ref{fig:SISyntheticViscosity}.

\subsection{Control of the heat equation}
\begin{figure}
  \centering
  \includegraphics[width=\linewidth]{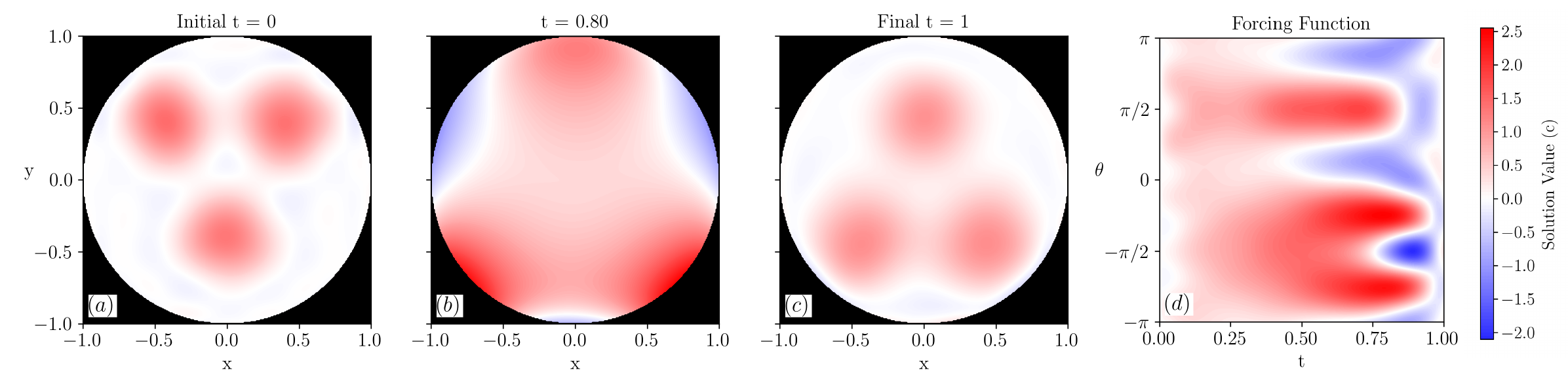}
  \caption{Heat equation solved on a circular embedded domain using $N=10$ basis functions in $x$, $y$, $t$, and $\theta$, where $\theta$ parameterizes position on the boundary. The heat equation was solved for $\alpha = 0.1$ where the value of the function on the boundary was prescribed by $u = f(\theta,t)$. The optimizer found a forcing function $f$ that started from an array of three dots (a) and inverted those dots by the end of the solve (c), with the optimized forcing function shown in (d).}
  \label{fig:heat_forcing}
\end{figure}
Besides inverting unknown fields from reference data, another important class of inverse problems is finding a control scheme to produce a target outcome. Our global solution representation lets us treat
a control input as an unknown field parameterized with its own set of basis functions and coefficients. We can then solve for this control function while {\sl simultaneously} solving the PDE. To do so, we augment our loss function in Eqn.~(\ref{eq:PINNLoss}) with targets for the underlying field at different time points. In practice, these targets are indistinguishable from a data residual, which both require the function to have some specified value at a particular coordinate.

To demonstrate the learning of a proper control sequence, we first consider the heat equation on the unit disk, embedded in $[-1,1]\times[-1,1]$. The governing equations are
\begin{gather}
	\frac{\partial u}{\partial t} = 0.1 \left( \frac{\partial^2 u}{\partial x^2} + \frac{\partial^2 u}{\partial y^2} \right),
\end{gather}
with initial conditions $u(x, y, t) = f(\theta,t)$ with $\theta = \tan^{-1}(y/x)$ on $x^2 +y^2 =1$. Here,
$f(\theta,t)$ is an unknown function of time defined on the boundary (parameterized by the angle $\theta$) that serves as the control input for this problem. 

We initialize the heat field in a pattern of three smooth dots arranged in a triangular pattern (Fig.~\ref{fig:heat_forcing}a), and attempt to find $f$ that causes the triangular arrangement to invert at $t=2$. We represent the heat field $u$ with a $N=10$ Chebyshev modes in each of $x$, $y$, and $t$, and express
the boundary function with with a $N=10$ mode Fourier basis in $\theta$, and $N=10$ Chebyshev modes in time. 

Fig.~\ref{fig:heat_forcing} shows the solution to this problem with (a), (b), and (c) show snapshots of the heat field at $t=0$, $t=1$, and $t=2$ respectively. The optimizer finds $f$ (Fig.~\ref{fig:heat_forcing}(d)) achieving the desired objective. The forcing function shows the emergence of the three-fold rotational symmetry, and note that the solution defines one lobe later in time due to the relatively closer position of the target triangle to the boundary along that component of the edge. 
What is striking about this particular example is that we are able to solve for both the control variable and the accompanying solution {\sl in a single optimization loop}. This contrasts with other differentiable PDE methods \cite{Alhashim2025}, where for each step in optimizing the control variable requires a complete run of the PDE solver to find a valid solution.

\subsection{Optimal transport}
\begin{figure}
  \centering
  \includegraphics[width=\linewidth]{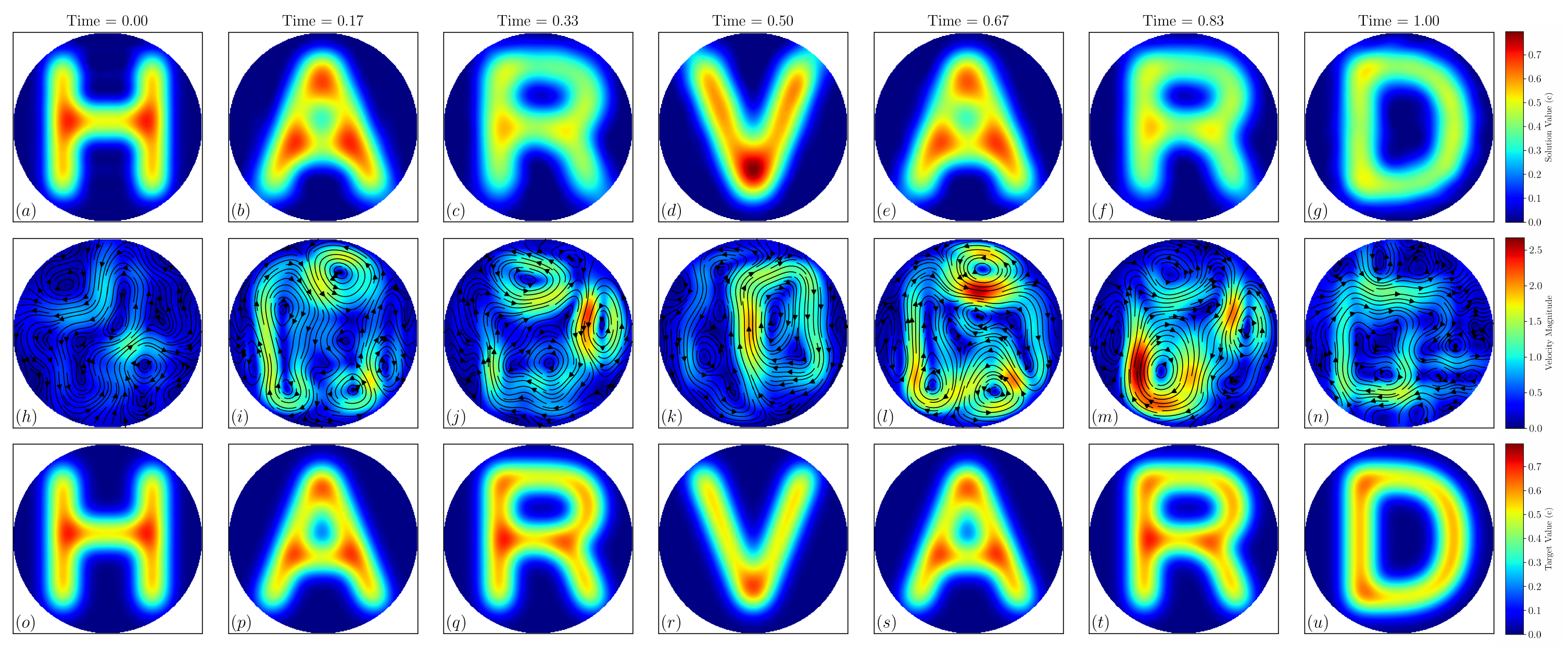}
  \caption{Top row: values of concentration field at different time points. Middle row: velocity field at different time points that transports concentration field into target shapes. Bottom row: target concentration field shapes specified at each discrete time point. Simulations were performed with a diffusion rate $D = 0.1$ with $N=15$ modes in each coordinate $x$, $y$, and $t$. Chebyshev polynomials were used to express the field in each direction. Solver conserved mass within 4\% across entire time range despite no explicit mass conservation residual term.}
  \label{fig:transport}
\end{figure}
As another example of solving an inverse control problem, we consider the advection-diffusion equations on the unit disk bounded by the unit circle $\partial \Omega$, 
\begin{gather}
	\frac{\partial c}{\partial t} + \vec{u} \cdot \nabla c = D \nabla^2 c,
\end{gather}
where $\vec{u}=(u,v)^\T$ is an incompressible velocity field ($\nabla \cdot \vec{u} = 0$) obeying no-slip and no-penetration conditions on the boundary, and the concentration $c$ has a no flux condition ($\frac{\partial c}{\partial n} = 0$) at the boundary. Our goal is to find a velocity field $\vec{u}$ that can transport the concentration field into a series of target shapes despite the influence of diffusion.

We solve this problem using $N=15$ Chebyshev modes in each of the three coordinates for $D = 0.1$ and $t\in[0,2]$. Our target profiles consist of seven individual targets at equally spaced time intervals, which together spell out Harvard, as shown in the bottom row of Fig.~\ref{fig:transport}. The individual letters were weighted to ensure they had the same total mass and were blurred to avoid sharp edges.

The solution concentration field and the optimized velocity field are shown in the top and middle row of Fig.~\ref{fig:transport} respectively. The optimizer finds a velocity function that leads to concentration profiles in good agreement with their targets.

Given the no-flux boundary conditions on $c$, the total mass $\int c \,\diff A$ should be conserved. When we evaluate the total mass across 200 evenly spaced time intervals we find that the solver conserved mass to within 4\% despite no explicit mass conservation residual. We expect this performance would be improved further with a higher number of basis modes and additional collocation points on the boundary.

\section{Discussion and Conclusion}

We have introduced a novel computational paradigm that marries the geometric flexibility of PINN-style global loss functions with the convergence properties of well-known spectral basis functions. By representing the solution as a coefficient tensor for a global basis, our approach solves forward, inverse, and control problems for PDEs on complex, irregular domains without the need for meshing. As demonstrated across a range of systems, from stiff, time-dependent systems to data assimilation on geophysical scales, this method can accurately find solutions to forward and inverse PDE problems.

This methodology offers significant gains in computational efficiency relative to other approaches for solving PDE-governed inverse problems. Our method requires orders of magnitude fewer parameters than comparable neural network approaches, enabling solutions to be found in minutes on consumer-grade hardware. This efficiency and accessibility significantly lower the barrier to data assimilation and inverse design. While not intended to replace highly optimized solvers for standard forward problems, which remain superior in both speed and accuracy, our approach excels in its flexibility for tackling a wide range of inverse, data assimilation, and control problems, particularly on complex geometries. Our approach also offers a route towards inverse design on complex problems \textit{without} having to translate solvers into differentiable programming languages (e.g. \cite{Alhashim2025}).
There are a few technical issues that need to be better understood with this approach.

First, while our method demonstrates empirical exponential convergence as the number of basis modes increases, solving for the coefficients of these bases using a global loss function defined over collocation points trades theoretical guarantees on finding the best coefficient representations for flexibility. 
Given that the minimization problem is in general a nonconvex function of the $\tensor{C}$, there are no guarantees for finding optimally convergent solution. The situation here is analogous to training deep neural networks, where there has been much work (e.g. \cite{kawaguchi2016deep,li2018visualizing}) understanding the role of local minima in loss landscapes. There, it has been noted that multiple local minima can be good solutions for the problem, even when evaluated on a held out test distribution. Further, while the convergence is exponential, it is too slow to be truly spectral, leaving the error relatively high compared to standard computational methods. We believe that this may be improved with a better choice of optimizer and preconditioner.

Analogously, while carrying out the work in this paper, we have found cases where there are multiple solutions minimizing the loss Eqn.~(\ref{eq:PINNLoss}) that exhibit exponential convergence. As an example, we consider our first problem, the Laplace problem on peanut geometry (Fig.~\ref{fig:laplace}). When we change the optimizer from a Dogleg algorithm to a regularized adaptive loss scheme with the ADAM optimizer \cite{Kingma2014}, we again observe exponential convergence (Fig.~\ref{fig:laplace_analysis}) but with a {\sl different} convergence rate. Indeed, the coefficient matrix found by each method is different with the ADAM solution representing a better (lower error) solution for a given representation.  This is striking, given that in the $N\to\infty$ limit, both methods converge to the same solution.

A careful examination of the optimized coefficient tensors reveals that many of the parameters often are unimportant for representing the solution. Fig.~\ref{fig:laplaceSVD} compares the singular values of the coefficient tensors of best Dogleg and best ADAM solutions from Fig.~\ref{fig:laplace_analysis}. First, it is clear that the ADAM solution represents the function in terms of much smaller coefficients, which is a consequence of the $L_1$ regularization term used in the loss function. However, the Dogleg solution shows a much steeper drop-off in magnitude of the singular values after the thirteenth, while the ADAM solution has a slower decay of the parameter values. The Dogleg singular values eventually plateau around the limit of single point machine precision.
Analysis of both of these sets of singular values indicates that in the case of Laplace's equation in the peanut a more efficient factorized version of the coefficient tensor could capture most of the variability of the solution. For reference, we have provided complete plots of the optimized coefficient tensors for every example shown in this paper in Figs.~\ref{fig:SILaplaceCoeffs} through \ref{fig:SIViscCoeffsB}.

In conclusion, by creating a single mesh-free and differentiable framework for solving PDE-constrained inverse problems, this work paves the way for tackling complex, multi-physics optimization and control problems that have long been out of reach for traditional solvers. While there is much work to be done, it seems reasonable to wonder whether these methods could open the door for an entirely different class of solvers for finding optimal solutions of PDEs to emerge.

\section{Methods}
We implement our method in Python using JAX \cite{jax2018github} and Equinox \cite{kidger2021equinox}, which enables both automatic differentiation and just-in-time compilation. To perform our minimization we use Optax \cite{deepmind2020jax} for first-order gradient-based methods and Optimistix \cite{optimistix2024} for second order Gauss-Newton type methods. Both packages make use of the automatic differentiation capabilities of JAX to perform efficient gradient computation.

\subsection{Loss function}
Our method finds an optimal set of coefficients $\tensor{C}$ by minimizing a composite loss function, Eqn.~(\ref{eq:PINNLoss}). Each individual term in the loss function is a mean over the residuals in each term. For example, consider a PDE system with a differential operator $\mathcal{D}_i$ with $i=1,...,k$ components. These components may be coupled systems of equations or different components of the operator acting on a vector of functions $\vec{u}=(u_1,\dots,u_n)^\T$. Given a set of $N_d$ collocation points, $\{\vec{x}_j\}_{j=1}^{N_d}$, the PDE loss term is given by
\begin{equation}
  \L_{\mathrm{pde}} = \frac{1}{kN_d}\sum_{i=1}^k\sum_{j=1}^{N_d}\mathcal{D}_i[\vec{u}(\vec{x_j})]^2.
\end{equation}
This averaging ensures that the loss is normalized with respect to the number of both operator components and collocation points. The boundary and data loss terms are computed similarly over their respective collocation points.

\subsection{Loss weighting}
In general, each component of the total loss may be weighted independently,
\begin{equation}
  \L = \lambda_1 \L_{\mathrm{pde}} + \lambda_2 \L_{\mathrm{boundary}} + \lambda_3\L_{\mathrm{data}}.
  \label{eq:PINNLossWeighted}
\end{equation}
Manually tuning the weights may be difficult, as the magnitude of each term can fluctuate significantly the optimization. We have found it beneficial in many cases to adaptively re-weight the terms, with the weights given by
\begin{equation}
  \lambda_i = \frac{\frac{1}{N} \sum_{j=1}^{N} ||\nabla L_j||}{||\nabla L_i|| + \epsilon}.
  \label{eq:reweighting}
\end{equation}
Here, $N$ is the total number of terms in the loss function, the gradients are with respect to the coefficient tensor $\tensor{
C}$, the $L_2$ norm is used, and $\epsilon$ is a small positive quantity to prevent numerical instability. This reweighting strategy ensures that loss terms with small residuals have equal weight in the overall gradient to terms with large residuals, keeping the different terms relatively balanced throughout the optimization. It also reduces the number of hyperparameters that need to be tuned to achieve a satisfactory solution.

\subsection{Optimization}\label{sec:optimizer}
For many problems second order Gauss-Newton methods were generally more reliable at converging to a solution without extensive hyperparameter tuning. In particular, we found success using both Dogleg and Levenberg–Marquardt algorithms on stiff problems and in the case of viscosity inversion. Dogleg tended to converge faster but was slightly less stable. We used the JAX compatible implementations of these two algorithms in Optimistix \cite{optimistix2024}, which would continue to step until a specified tolerance was reached. We set the relative and absolute tolerance in these problems to $10^{-6}$. When using second order methods we found that we could generally set all of the loss weights $\lambda_i = 1.0$.

In problems involving higher numbers of dimensions, such as solving for the time evolution of some function over a two-dimensional surface, using second order methods can become impractical due to the computational overhead in finding estimates for the Hessian matrix. In these cases, we used first order JAX-compatible ADAM and NADAM algorithms as implemented in Optax \cite{deepmind2020jax}. 
In general, we found that these methods were much less stable at converging to a plausible solution than second order methods. To help address this, when using first order solvers we almost always used the adaptive weighting scheme described in the previous section. An $L_1$ penalty on the values of the coefficients themselves also improved optimizer stability and lead to sparser solutions.

We found that the optimizers were much more stable when we preconditioned the coefficient tensors. While there are likely problem-specific preconditioning that would lead to better results, we found that decent optimization performance could be achieved with a relatively simple preconditioner, with weights given by
\begin{equation}
  w(k_1,...,k_d) = \frac{1}{1+k_1^2+...+k_d^2},
\end{equation}
where $k_i$ is the mode index for the basis functions in each direction.

In addition to the optimizer parameters, an important consideration is the number of collocation points chosen for a particular problem. In general, higher numbers of modes require more collocation points to find a satisfactory fit, which in turns leads to higher memory overheads in the optimization. We sampled collocation points once before beginning the optimization and did not do any resampling, though this may improve performance. We sampled collocation points for equation, boundary, and data randomly. In cases of complex boundaries defined by a mask from data, we sampled a subset of the grid points defined by the data. In cases where the domain was defined analytically we uniformly sampled over the entire domain.

For a complete summary of the optimizer types and hyperparameters used in producing each example seen in the text, including the number of collocation points, refer to Table \ref{tab:my_new_table}. For a comparison of how well different combinations of solvers and hyperparameters are able to converge as the number of modes $N$ increases, refer to Fig.~\ref{fig:laplace_analysis}.
\subsection{Laplace Equation}
The precise form of the boundary in Fig.~\ref{fig:laplace} is given by
\begin{equation}
  r(\theta)=\frac{0.45 \cos(2\theta) + 0.55}{(0.375 \cos(2\theta) + 0.625)^{\frac{3}{2}}},
  \label{eq:peanut}
\end{equation}
where $\theta = \tan^{-1}(y/x)$ is the standard polar angle. 

Our domain also includes a circular cut out of radius $0.15$ at the location $(.3,.1)$ which defines the inner boundary $\partial \Omega_2$. We set $u$ on the outer boundary as $u(\partial\Omega_1) = 2\sin(\theta) + \cos(3\theta)$, while
on the inner boundary $u(\partial\Omega_2) = -2$.

\section{Acknowledgements}
We thank Francesco Mottes for important discussions and Nick Trefethen for critical feedback. This work was supported by the Office of Naval Research through grant number ONR N00014-23-1-2654 and the NSF AI Institute of Dynamic Systems
(2112085).
\bibliography{arxiv_submission}

\clearpage
\onecolumngrid
\appendix
\setcounter{section}{0}
\setcounter{figure}{0}
\setcounter{table}{0}
\setcounter{equation}{0}

\renewcommand{\thesection}{S\arabic{section}}
\renewcommand{\thefigure}{S\arabic{figure}}
\renewcommand{\thetable}{S\arabic{table}}
\renewcommand{\theequation}{S\arabic{equation}}

% If you use hyperref, also reset the PDF anchors so links are correct:
\makeatletter
\renewcommand{\theHfigure}{S\arabic{figure}}
\renewcommand{\theHtable}{S\arabic{table}}
\renewcommand{\theHequation}{S\arabic{equation}}
\makeatother
\section{Supplemental Information}
\subsection{Validation of data-assimilation for ice sheets viscosity inversion}
We present a solution for the inversion of a viscosity field on an ice sheet through the shallow shelf approximation (SSA) equations based on experimental measurements of the thickness $h$ and velocity components $u$ and $v$. The SSA equations are given in terms of these variables and a spatially varying viscosity $\mu$ in Eqns.~(\ref{eq:ssa_eq1}) and (\ref{eq:ssa_eq_2}), which we reproduce here as
\begin{equation}
\frac{\partial}{\partial x}\left(4\mu h \frac{\partial u}{\partial x} + 2\mu h \frac{\partial v}{\partial y}\right) + \frac{\partial}{\partial y}\left(\mu h \frac{\partial u}{\partial y} + \mu h \frac{\partial v}{\partial x}\right) \\= \rho_i g \left(1 - \frac{\rho_i}{\rho_w}\right) h \frac{\partial h}{\partial x},
\end{equation}
\begin{equation}
\frac{\partial}{\partial y}\left(4\mu h \frac{\partial v}{\partial y} + 2\mu h \frac{\partial u}{\partial x}\right) + \frac{\partial}{\partial x}\left(\mu h \frac{\partial u}{\partial y} + \mu h \frac{\partial v}{\partial x}\right) = \rho_i g \left(1 - \frac{\rho_i}{\rho_w}\right) h \frac{\partial h}{\partial y}.
\end{equation}
$\rho_i$ and $\rho_w$ are the density of ice and water, taken to be 917 and 1030 $\mathrm{kg}/\mathrm{m}^3$. In addition to these equations, there is a force balance condition (Eqns.~(\ref{eq:ssa_bc_1}) and (\ref{eq:ssa_bc_2})) applied at the calving front, 
\begin{equation}
2\mu h \left(2\frac{\partial u}{\partial x} + \frac{\partial v}{\partial y}\right) n_x + \mu h \left(\frac{\partial u}{\partial y} + \frac{\partial v}{\partial x}\right) n_y = \frac{1}{2}\rho_i g \left(1 - \frac{\rho_i}{\rho_w}\right) h^2 n_x,
\end{equation}
\begin{equation}
\mu h \left(\frac{\partial u}{\partial y} + \frac{\partial v}{\partial x}\right) n_x + 2\mu h \left(\frac{\partial u}{\partial x} + 2\frac{\partial v}{\partial y}\right) n_y = \frac{1}{2}\rho_i g \left(1 - \frac{\rho_i}{\rho_w}\right) h^2 n_y.
\end{equation}
The calving front has a normal vector $\vec{n}=(n_x,n_y)^T$. 

In the main text, we show our solution applied to the Amery ice sheet using data from \cite{Wang2025}. In addition to the real ice sheet data, the authors also provide synthetic data generated in COMSOL for a solution to these equations to validate their methodology for viscosity inversion, given that the viscosity field of the COMSOL data is prescribed. 

\begin{figure}[h!]
  \centering
  \includegraphics[width=\linewidth]{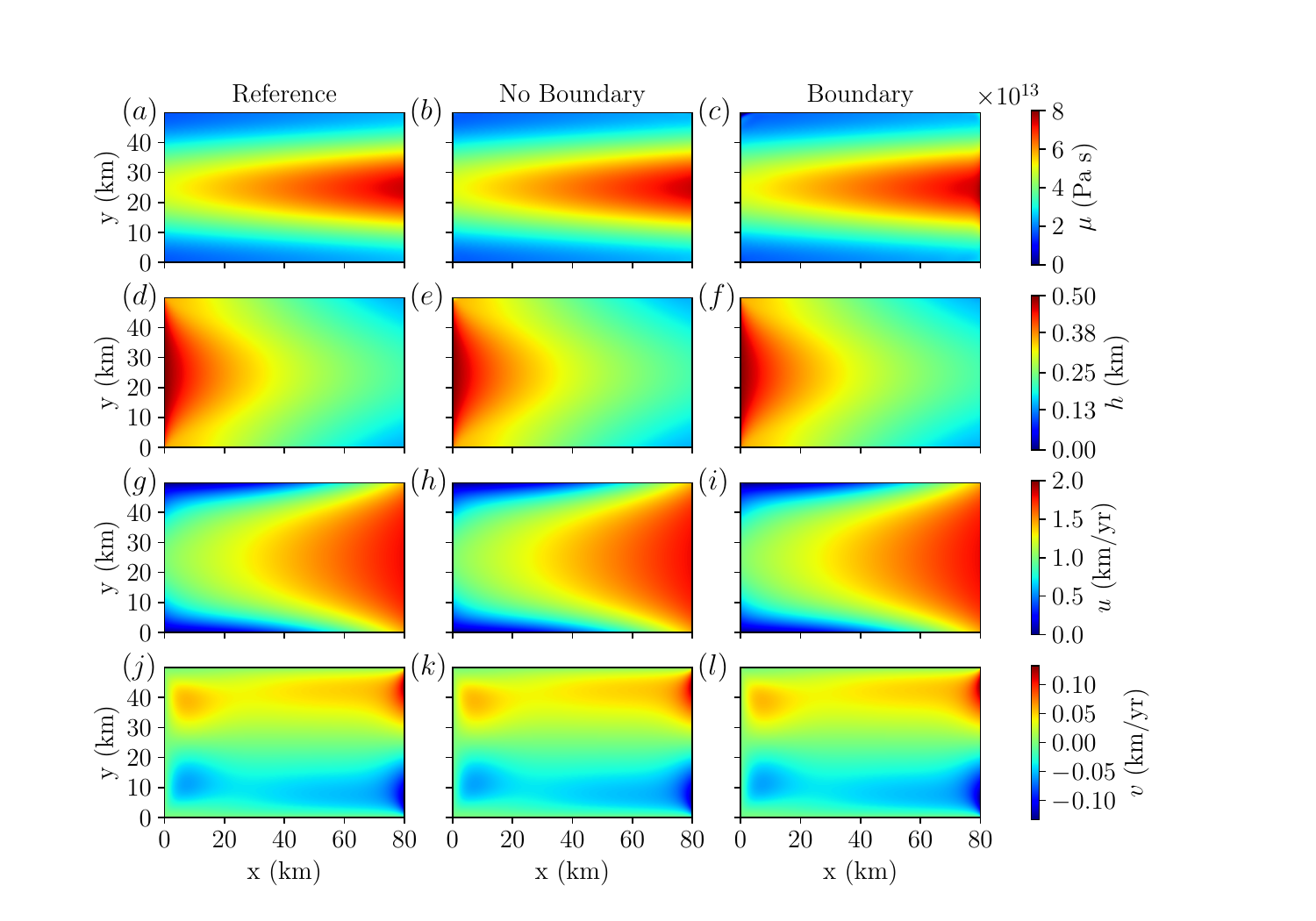}
  \caption{Validation of viscosity inversion for ice sheets governed by the shallow shelf approximation. The first column shows the reference COMSOL solution from \cite{Wang2025}. The other columns show the results of solving the SSA equations with access to reference data in $h$, $u$, and $v$ but with no data on the viscosity $\mu$, which must be inverted. The second column shows the results of solving this inversion problem without the calving front boundary condition while the third column shows the results with this boundary condition.}
  \label{fig:SISyntheticViscosity}
\end{figure}

We apply our methodology to invert the viscosity from the synthetic data, using data for $h$, $u$, and $v$ as a reference as in the full ice sheet problem. As before, we fit coefficient matrices for $\mu$ as well as $h$, $u$, and $v$ that ensure the equations are satisfied while interpolating the provided data. We solve the system both with and without the boundary condition, with results shown in Fig.~\ref{fig:SISyntheticViscosity}. We use Chebyshev polynomials to represent the solution, with $N=20$ basis modes in each direction. We solve using a Levenberg-Marquardt algorithm, which is a second-order Gauss Newton-type solver like was used in generating Fig.~\ref{fig:ameryinversion}. We used 2,000 PDE and data points and 600 points on the boundary. Each term in the loss was assigned a weight of one.

We are able to accurately invert the viscosity field, with an $L_2$ error calculated on the non-dimensionalized viscosity of $0.005$ when the boundary condition is included and $0.0008$ when it is not. This is in contrast to the PINN validation shown in \cite{Wang2025}, which was unable to solve for the viscosity field without the boundary condition. The increased representational efficiency of our method means that the data itself is sufficient to ensure inversion of the true viscosity field.
\clearpage
\section{Supplemental Figures}
\begin{figure}
  \centering
  \includegraphics[width=\linewidth]{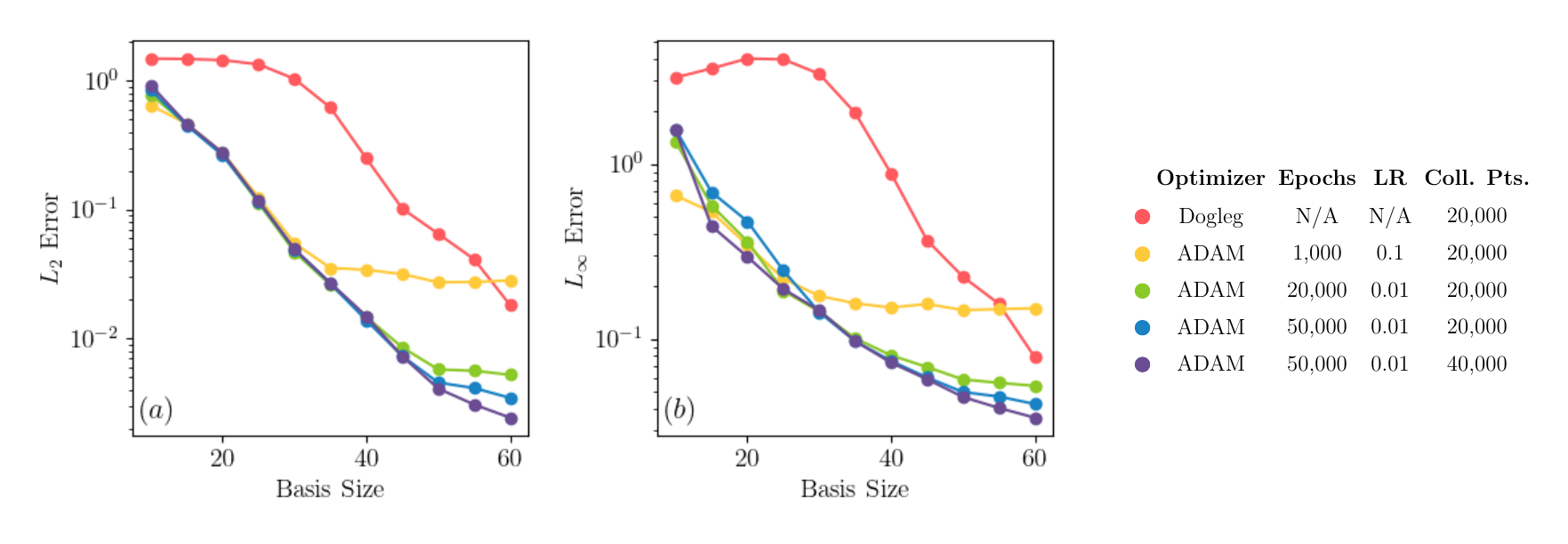}
  \caption{$L_2$ (a) and $L_\infty$ errors for solutions to Laplace's equations solved on the peanut shape. The result shown in Fig.~\ref{fig:laplace}(c) is in red. The other results present optimizations run with ADAM and a regularized adaptive loss function with various epochs, learning rates, and numbers of collocation points. Both methods show evidence of exponential convergence in the representation but each optimizer converges on different solutions. As the number of basis modes increases, more epochs and collocation points become necessary to continue the exponential convergence in the representation.}
  \label{fig:laplace_analysis}
\end{figure}
\begin{figure}
  \centering
  \includegraphics[width=0.5\linewidth]{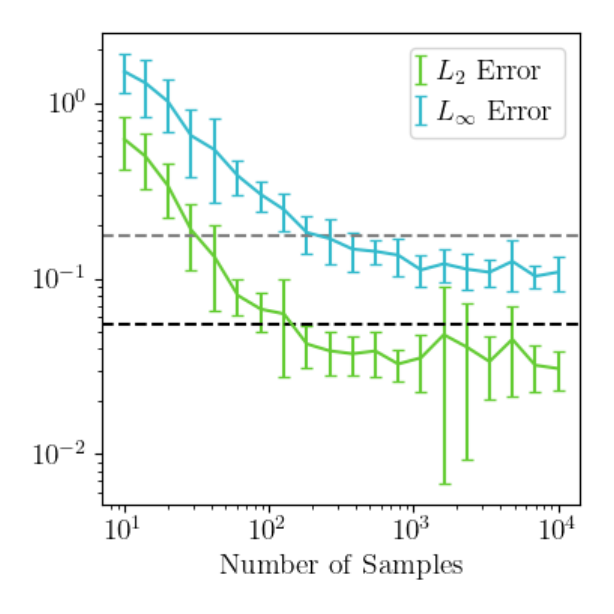}
  \caption{Mean and standard deviation of $L_2$ and $L_\infty$ (b) errors for solving Laplace's equation with $N=30$ Chebyshev polynomials in both of $x$ and $y$ on the peanut shape with a hole using a data fitting residual in place of a boundary condition residual. We repeatedly solved the equation using 15 repetitions of resampling the data for a given number of samples. In these examples, we use 20,000 collocation points, 1,000 epochs, and a learning rate of 0.1 with an ADAM optimizer, as in the yellow example from Fig.~\ref{fig:laplace_analysis}. The $L_2$ and $L_\infty$ error from Fig.~\ref{fig:laplace_analysis} for the equivalent solution is shown as the black and gray dashed lines respectively. We observe that the solution saturates at a particular error level as the number of samples increases but that this error is lower than when using the boundary residual.}
  \label{fig:laplaceinverseerror}
\end{figure}
\begin{figure}
  \centering
  \includegraphics[width=0.7\linewidth]{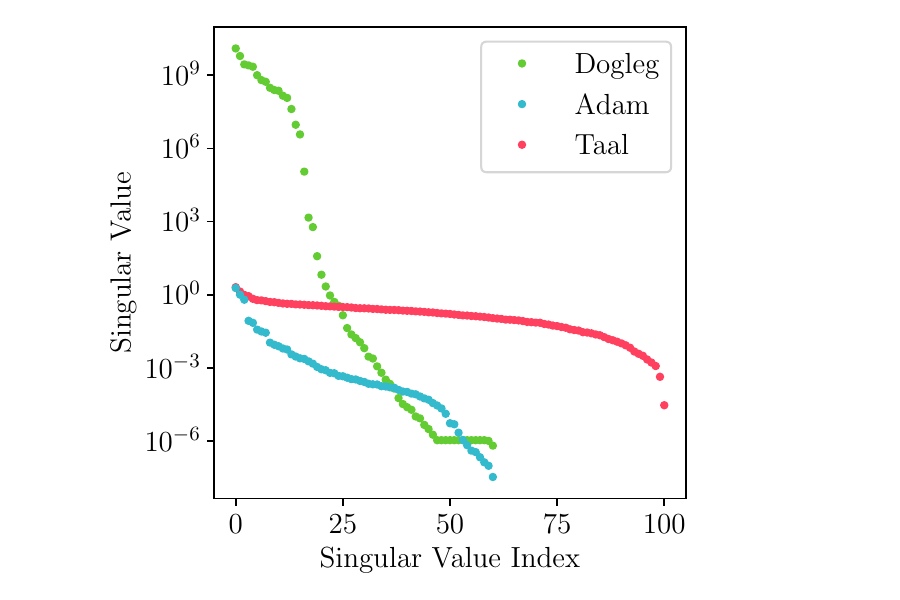}
  \caption{Singular values for the optimized coefficient matrices for solutions to Laplace's equation. The data labeled Dogleg and ADAM represent the two best solutions for their respective methods shown in Fig.~\ref{fig:laplace_analysis}, which both go up to 61 total values. Dogleg tends to push the solution to much higher singular values while having a much steeper drop off than ADAM, giving a more efficient representation (at the cost of accuracy). The dataset labeled Taal corresponds to the singular values of the solution to Laplace's equation on Lake Taal, shown in Fig.~\ref{fig:laplace}(f). The increased complexity from the boundary shape requires a much larger number of singular values represent the solution.}
  \label{fig:laplaceSVD}
\end{figure}
\begin{figure}
  \centering
  \includegraphics[width=\linewidth]{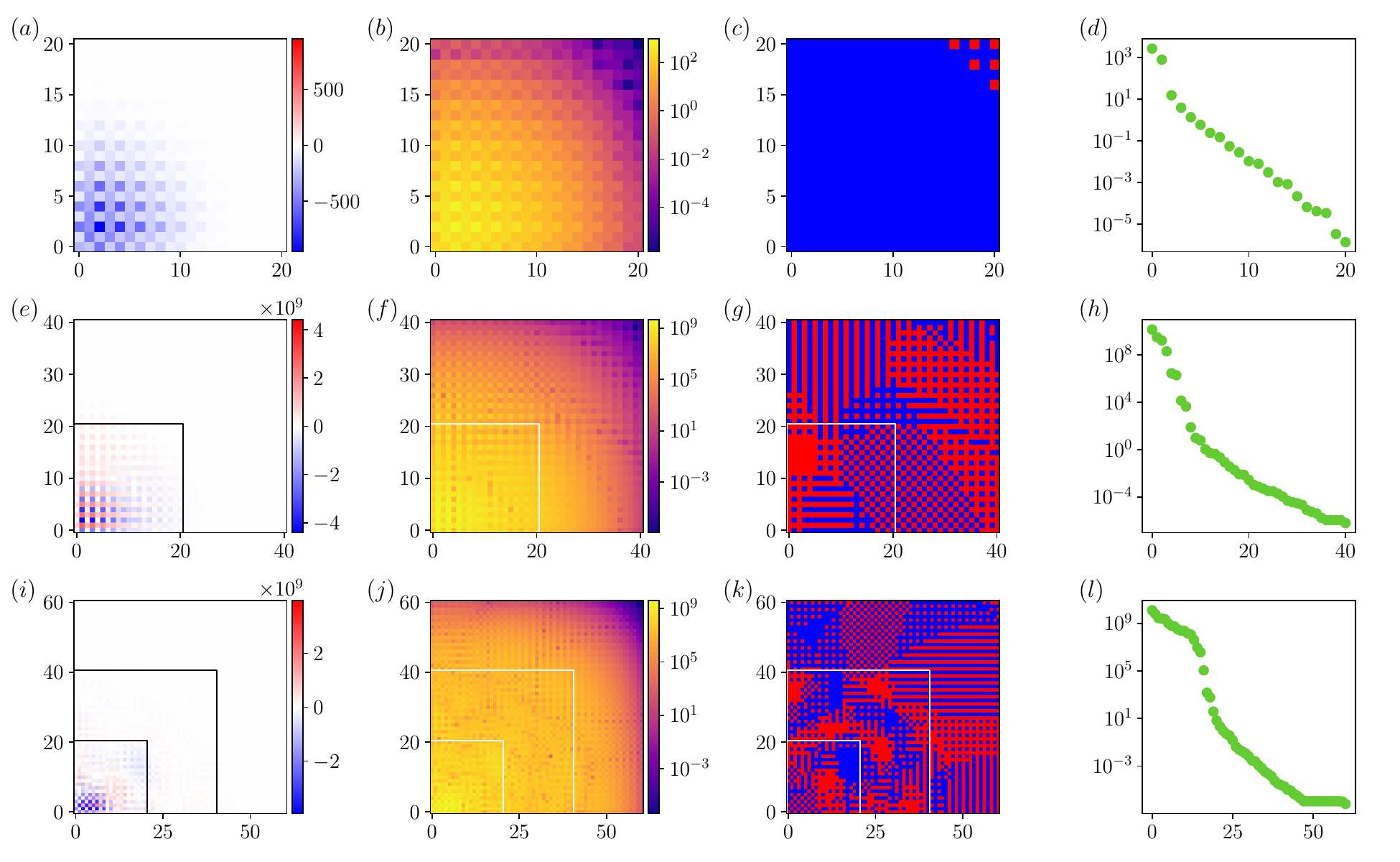}
  \caption{Coefficient matrices for the solution to Laplace's equation in the peanut shaped geometry shown in Fig.~\ref{fig:laplace}. These solutions were generated with the Dogleg solver. The first row shows the solution for 21 total basis modes, the second row for 41 total basis modes, and the final row for 61 total modes. The first column shows the raw coefficient values, the second column the absolute value of the coefficients plotted on a log scale, and the third column the sign of the coefficients. Red indicates positive values while blue represents negative ones. The fourth columns shows the singular values of the coefficient matrices.}
  \label{fig:SILaplaceCoeffs}
\end{figure}
\begin{figure}
  \centering
  \includegraphics[width=\linewidth]{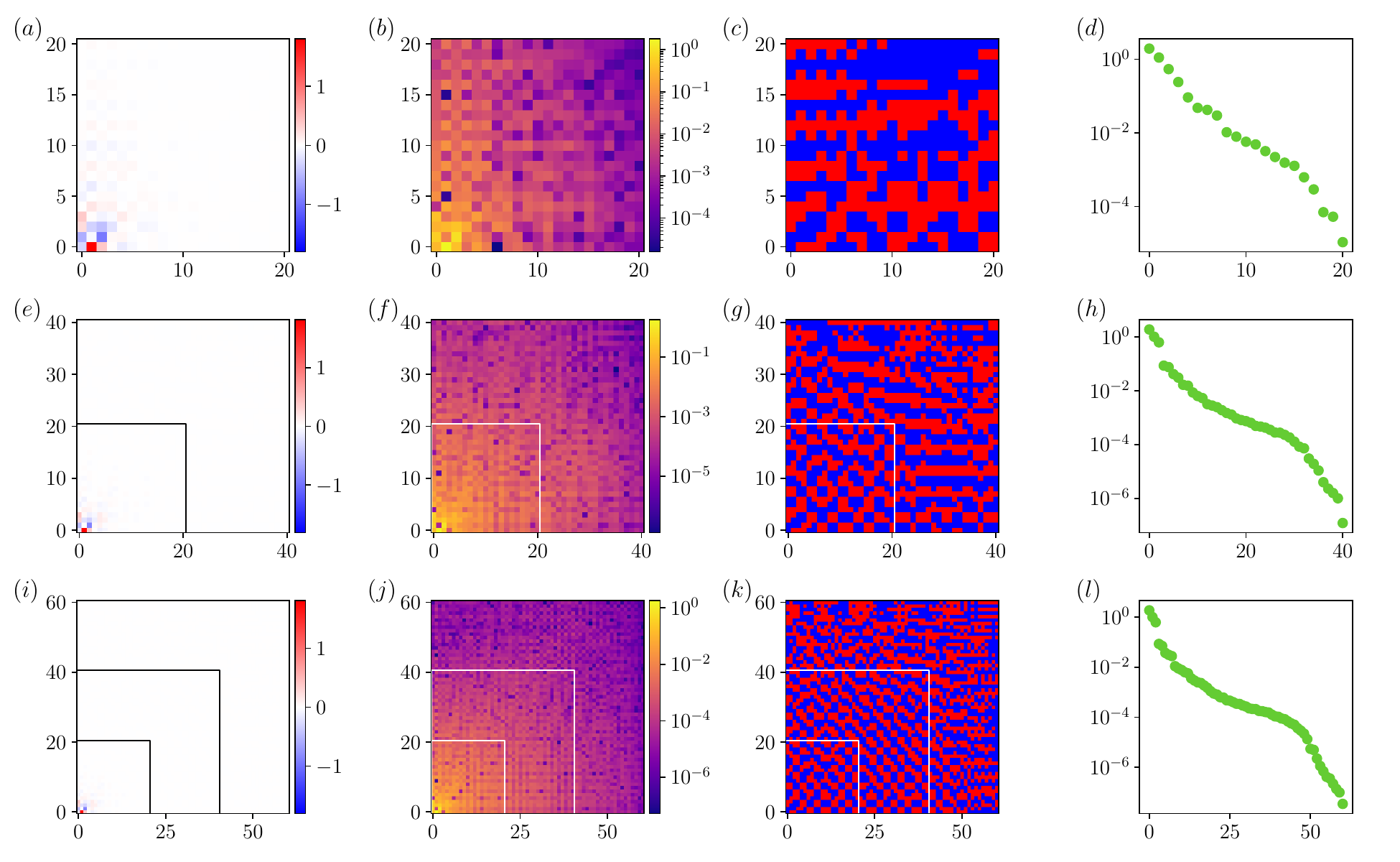}
  \caption{Coefficient matrices for the solution to Laplace's equation in the peanut shaped geometry shown in Fig.~\ref{fig:laplace}. These solutions were generated with the best ADAM solver shown in Fig.~\ref{fig:laplace_analysis}. The first row shows the solution for 21 total basis modes, the second row for 41 total basis modes, and the final row for 61 total modes. The first column shows the raw coefficient values, the second column the absolute value of the coefficients plotted on a log scale, and the third column the sign of the coefficients. Red indicates positive values while blue represents negative ones. The fourth columns shows the singular values of the coefficient matrices. The regularization term clearly drives the coefficients to much smaller values than in the Dogleg case.}
\label{fig:SILaplaceCoeffsAdam}
\end{figure}
\begin{figure}
  \centering
  \includegraphics[width=\linewidth]{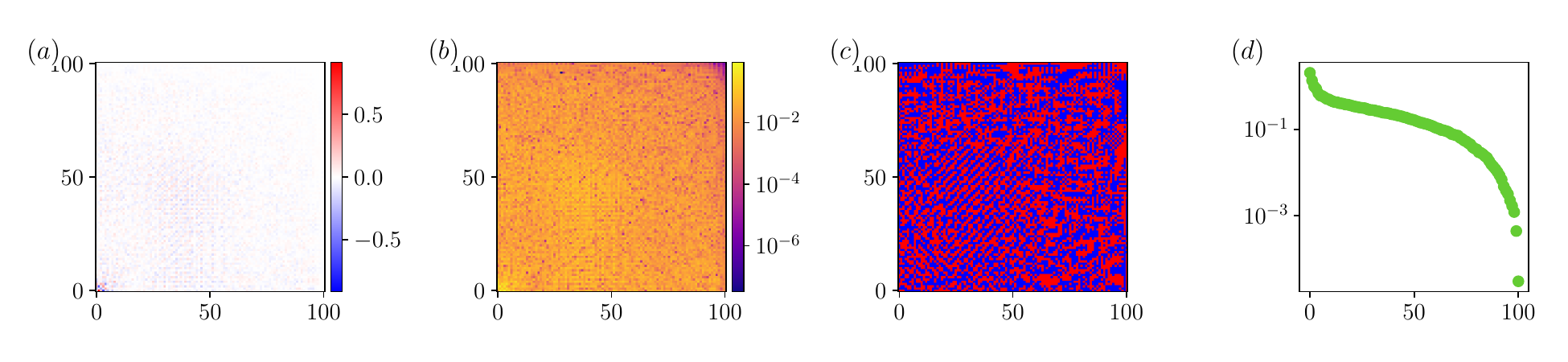}
  \caption{Coefficient matrices for the solution to Laplace's equation in Lake Taal as shown in Fig.~\ref{fig:laplace}(f). These solutions were generated with the Dogleg solver for 101 total basis modes. The first column shows the raw coefficient values, the second column the absolute value of the coefficients plotted on a log scale, and the third column the sign of the coefficients. Red indicates positive values while blue represents negative ones. The fourth columns shows the singular values of the coefficient matrices. When compared with the solutions for Laplace's equation in a peanut (Figs.~\ref{fig:SILaplaceCoeffs} and \ref{fig:SILaplaceCoeffsAdam}), it is clear that far more of the modes are important in capturing the complete solution.}
  \label{fig:SITaalCoeffs}
\end{figure}
\begin{figure}
  \centering
  \includegraphics[width=\linewidth]{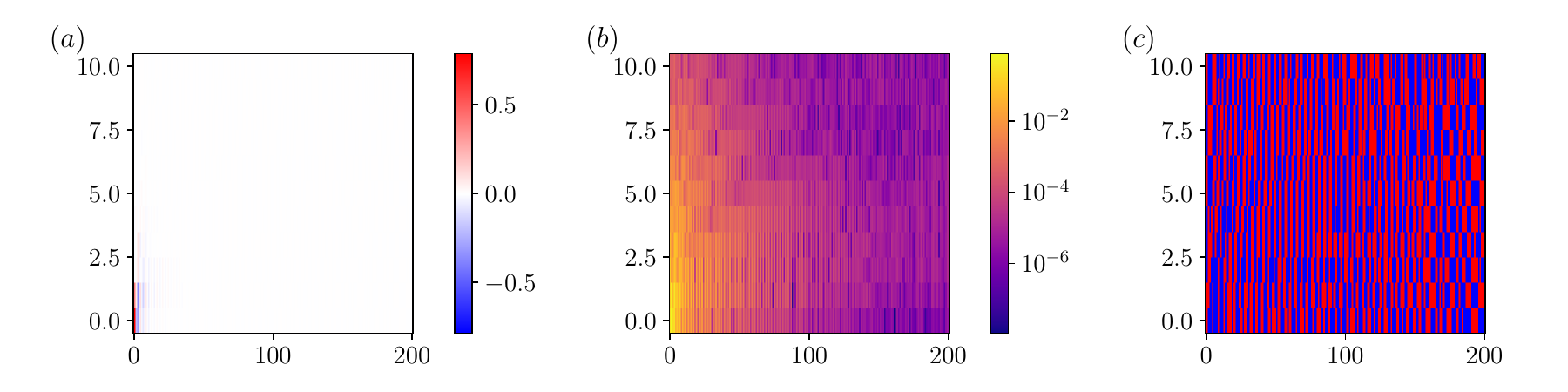}
  \caption{Coefficient matrices for the solution to the Allen-Cahn equation, as shown in Fig.~\ref{fig:stiff_eqs}(a). These solutions were generated with the Dogleg solver for 10 modes in time and 200 modes in the spatial variable $x$. The first column shows the raw coefficient values, the second column the absolute value of the coefficients plotted on a log scale, and the third column the sign of the coefficients. Red indicates positive values while blue represents negative ones.}
  \label{fig:SIACCoeffs}
\end{figure}
\begin{figure}
  \centering
  \includegraphics[width=\linewidth]{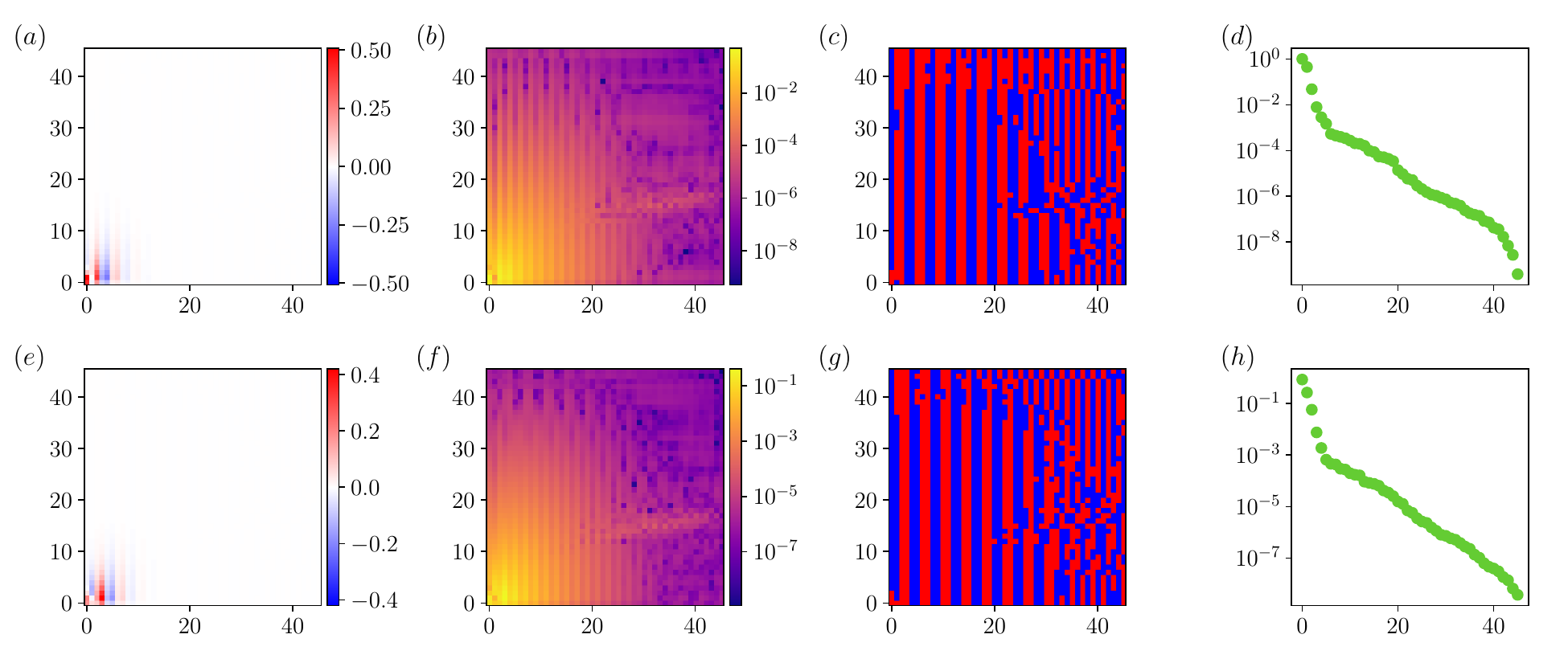}
  \caption{Coefficient matrices for the solution to nonlinear Schr\"odinger equation as shown in Fig.~\ref{fig:stiff_eqs}(d). These solutions were generated with the Dogleg solver with 46 total basis modes. The first row shows the solution for $u$ while the second gives the solution for $v$. The first column shows the raw coefficient values, the second column the absolute value of the coefficients plotted on a log scale, and the third column the sign of the coefficients. Red indicates positive values while blue represents negative ones. The fourth columns shows the singular values of the coefficient matrices.}
  \label{fig:SINLSCoeffs}
\end{figure}
\begin{figure}
  \centering
  \includegraphics[width=\linewidth]{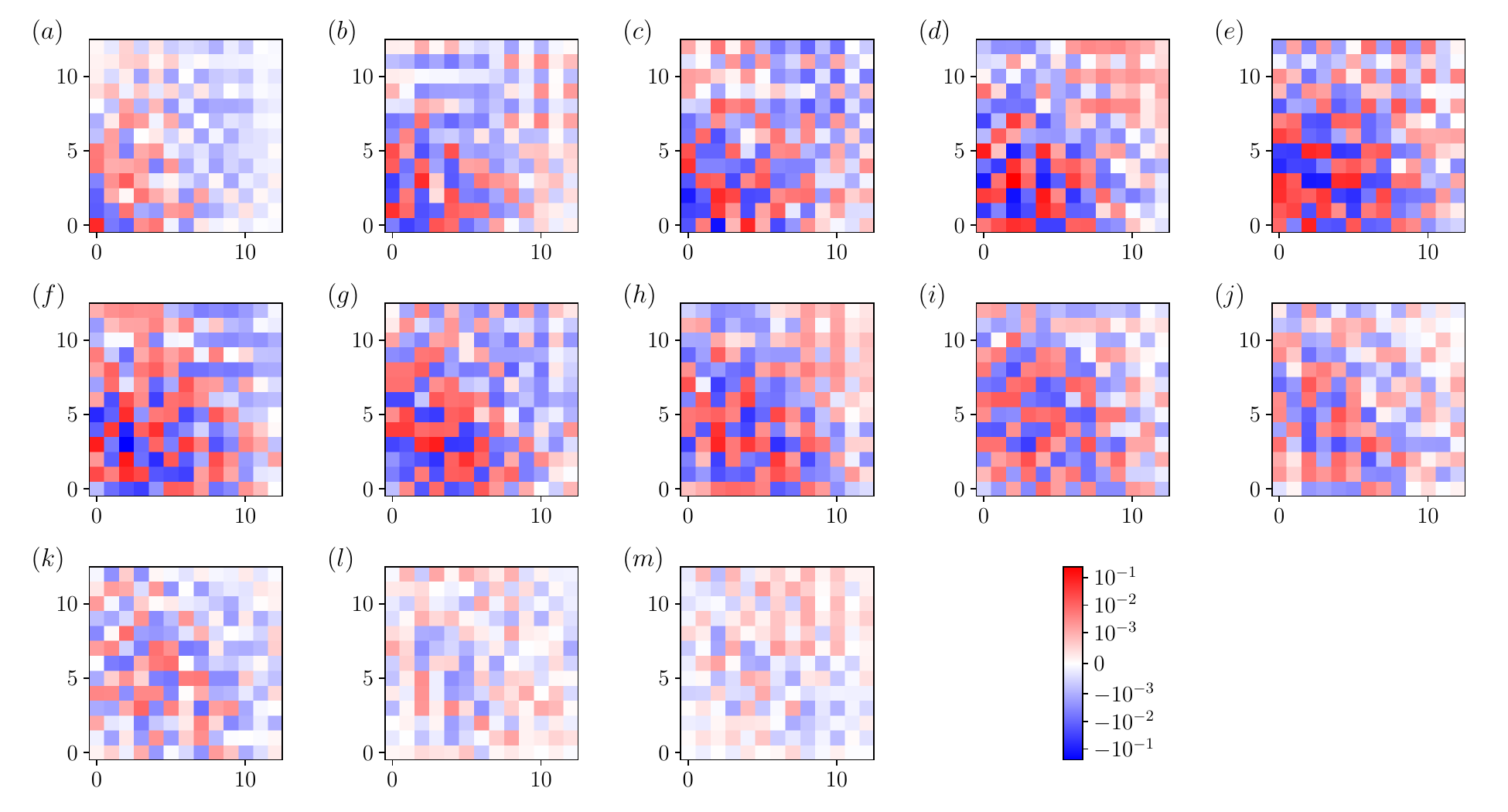}
  \caption{Coefficient matrix for solutions to the wave equation on the peanut geometry, as shown in Fig.~\ref{fig:wave}. The coefficients were found using an ADAM solver with a regularized adaptive loss. Here, each image is one slice of the total (13,13,13) tensor. The values are shown on a symmetric log scale.}
  \label{fig:SIWaveCoeffs}
\end{figure}
\begin{figure}
  \centering
  \includegraphics[width=\linewidth]{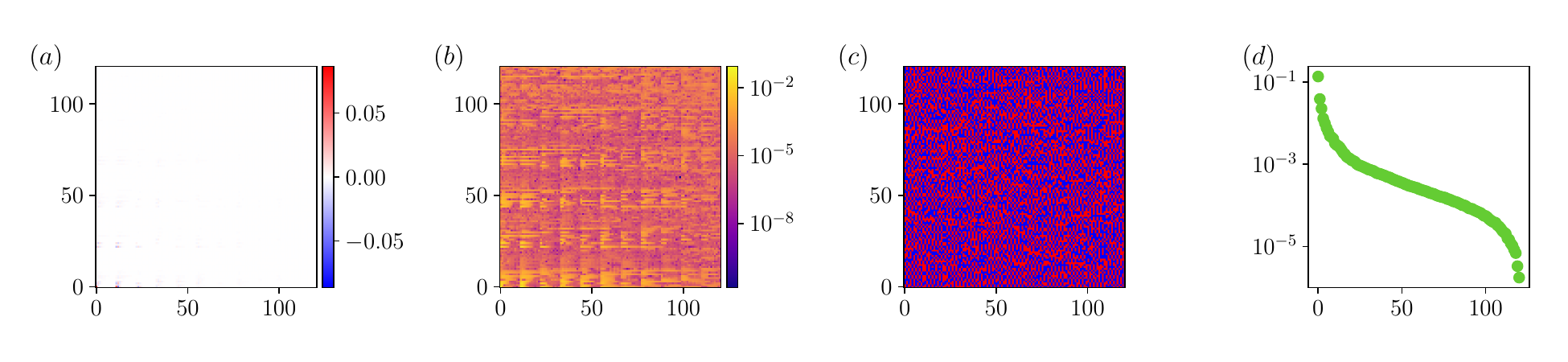}
  \caption{Coefficient matrix for the solution to the diffusion equation on the surface of a sphere, as shown in Fig.~\ref{fig:spherediff}. The coefficients were found using an ADAM solver with a regularized adaptive loss. The original tensor had dimensions (11,11,11,11) and has been flattened into (121,121) for display, leading to a checkerboard pattern in the figure. The first column shows the raw coefficient values, the second column the absolute value of the coefficients plotted on a log scale, and the third column the sign of the coefficients. Red indicates positive values while blue represents negative ones. The fourth columns shows the singular values of the coefficient matrices.}
  \label{fig:SISphereCoeffs}
\end{figure}
\begin{figure}
  \centering
  \includegraphics[width=\linewidth]{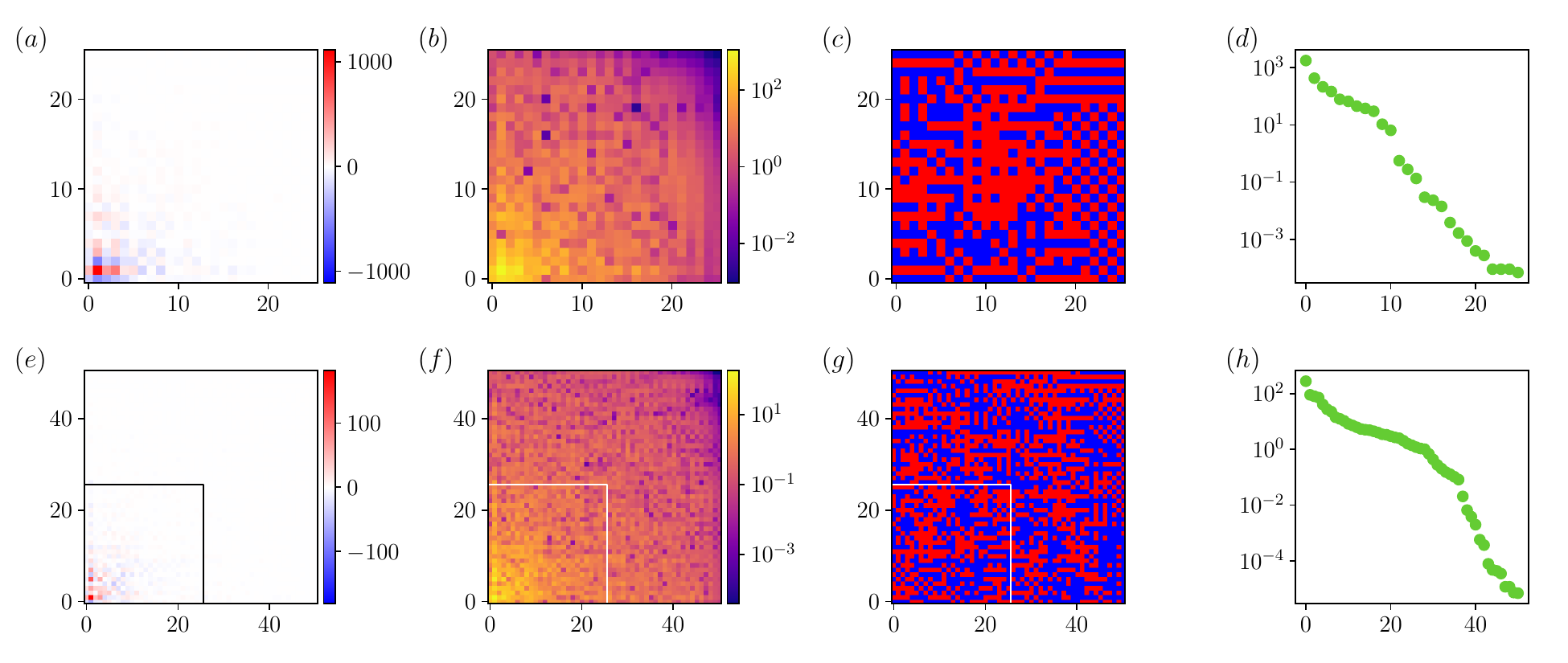}
  \caption{Coefficient matrices for the inverted viscosity $\mu$ from the shallow shelf approximation equations on the Amery ice sheet as shown in Fig.~\ref{fig:ameryinversion}. These solutions were generated with the Dogleg solver. The first row shows the solution for 26 total basis modes and the second row for 51 total basis modes. The first column shows the raw coefficient values, the second column the absolute value of the coefficients plotted on a log scale, and the third column the sign of the coefficients. Red indicates positive values while blue represents negative ones. The fourth columns shows the singular values of the coefficient matrices.}
  \label{fig:SIAmeryCoeffs}
\end{figure}
\begin{figure}
  \centering
  \includegraphics[width=\linewidth]{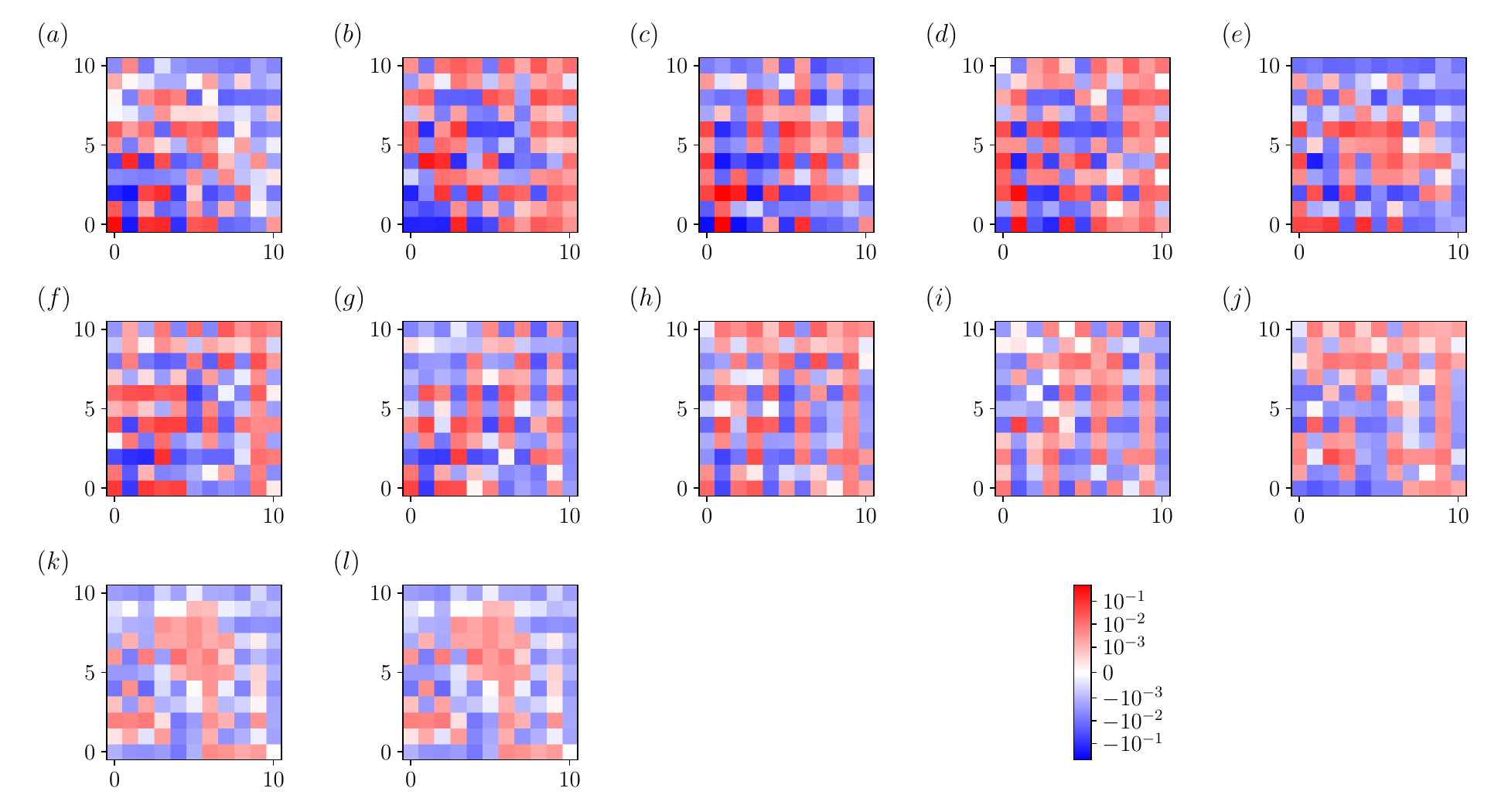}
  \caption{Coefficient matrix for solutions to the heat forcing problem, as shown in Fig.~\ref{fig:heat_forcing}. The coefficients were found using an ADAM solver with a regularized adaptive loss. Here, figures (a)-(k) is one slice of the total (11,11,11) tensor. Figure (l) shows the values of the coefficients for the boundary forcing function. The values are shown on a symmetric log scale.}
  \label{fig:SIHeatCoeffs}
\end{figure}
\begin{figure}
  \centering
  \includegraphics[width=\linewidth]{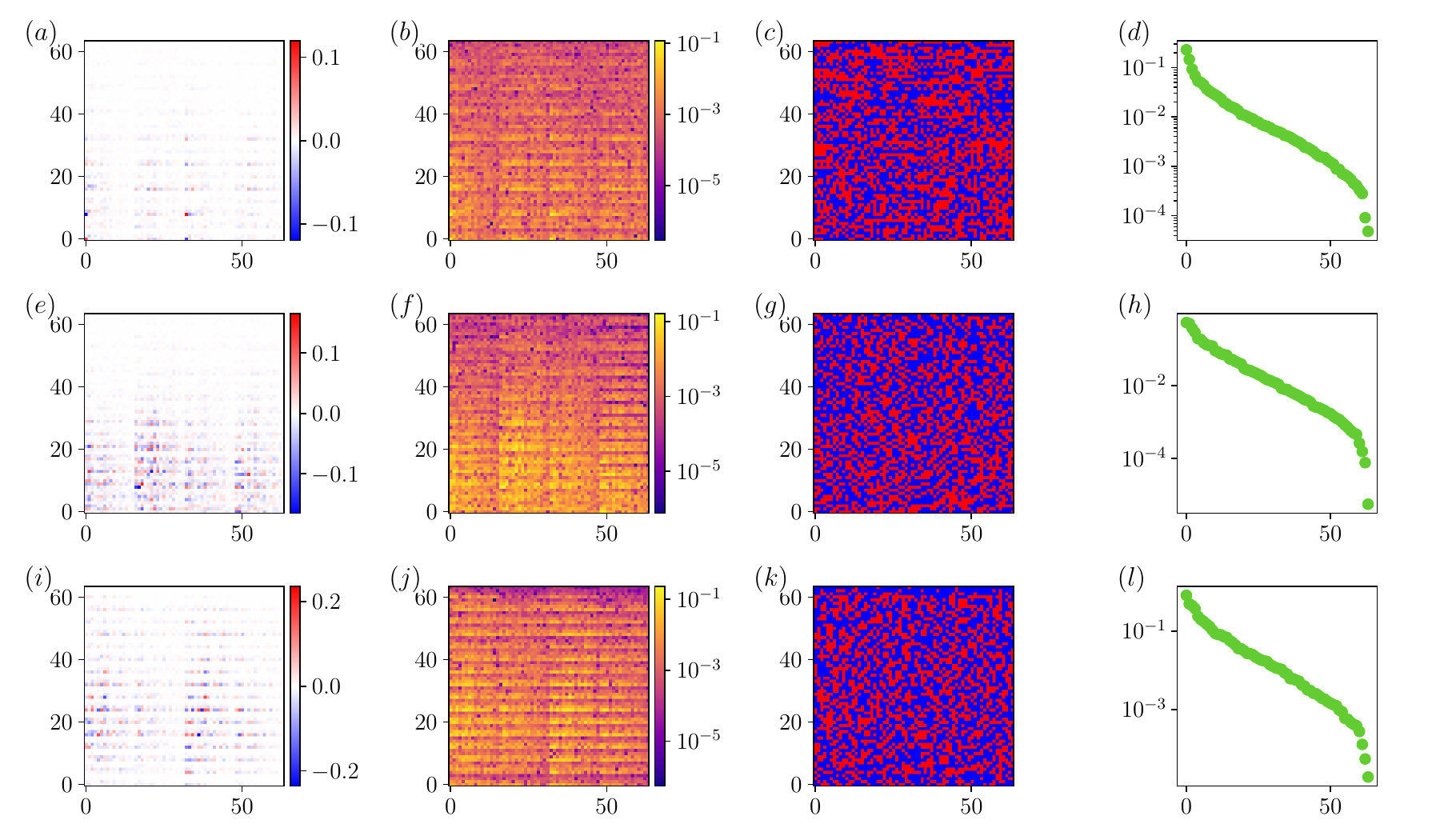}
  \caption{Coefficient matrix for the solution to the advection diffusion equation, as shown in Fig.~\ref{fig:transport}. The coefficients were found using an ADAM solver with a regularized adaptive loss. The original tensors had dimensions (16,16,16) and have been flattened into (64,64) for display, leading to a checkerboard pattern in the figure. The top row shows the coefficients for $c$, the middle for $u$, and the bottom for $v$. The first column shows the raw coefficient values, the second column the absolute value of the coefficients plotted on a log scale, and the third column the sign of the coefficients. Red indicates positive values while blue represents negative ones. The fourth columns shows the singular values of the coefficient matrices.}
  \label{fig:SITransportCoeffs}
\end{figure}
\begin{figure}
  \centering
  \includegraphics[width=\linewidth]{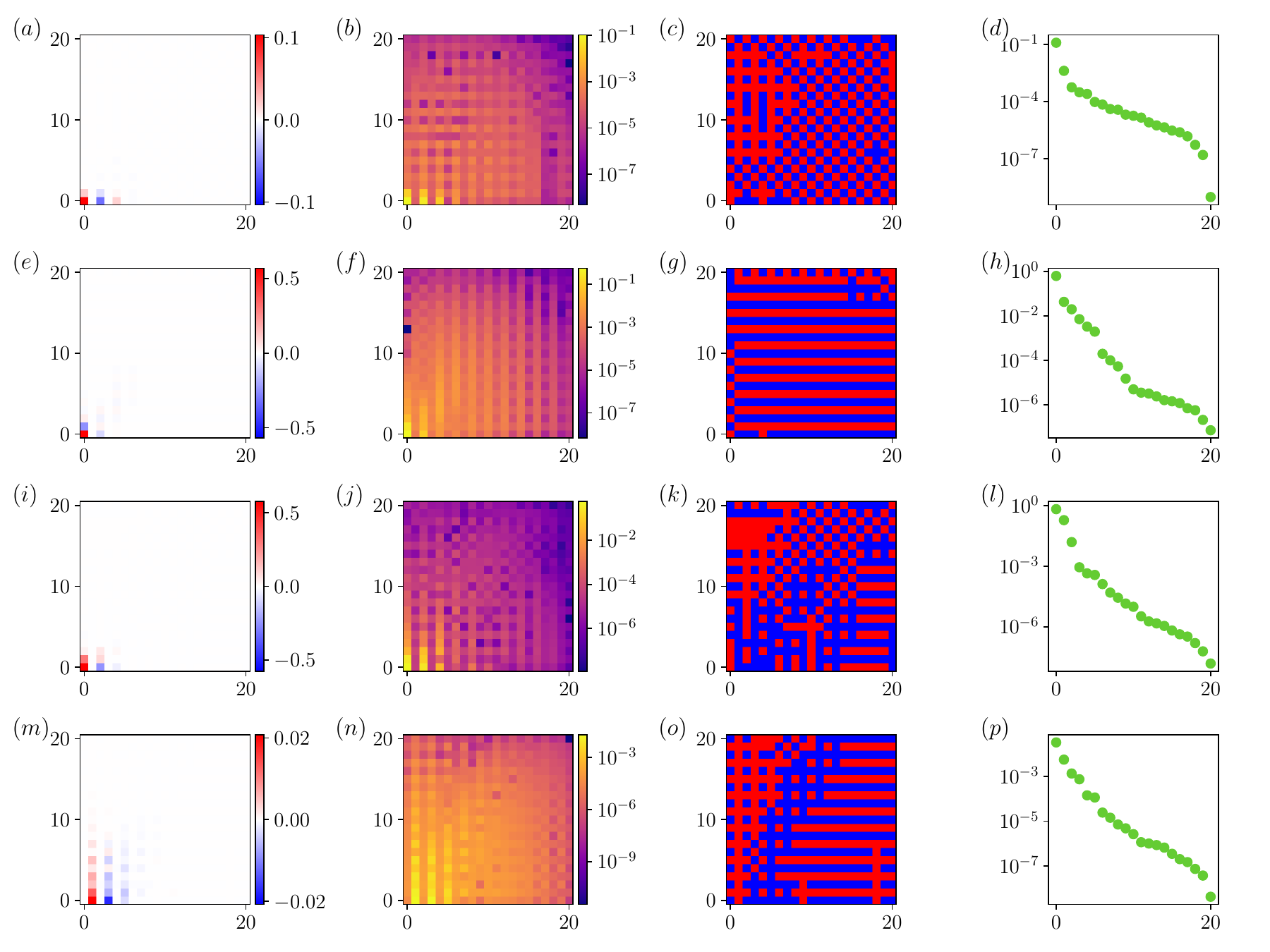}
  \caption{Coefficient matrix for the viscosity inversion solution with no boundary condition shown in the second column of Fig.~\ref{fig:SISyntheticViscosity}. The coefficients were found using a Levenberg Marquardt solver. The top row shows the coefficients for $\mu$, the second for $h$, the third for $u$, and the bottom for $v$. The first column shows the raw coefficient values, the second column the absolute value of the coefficients plotted on a log scale, and the third column the sign of the coefficients. Red indicates positive values while blue represents negative ones. The fourth columns shows the singular values of the coefficient matrices.}
  \label{fig:SIViscCoeffsNB}
\end{figure}
\begin{figure}
  \centering
  \includegraphics[width=\linewidth]{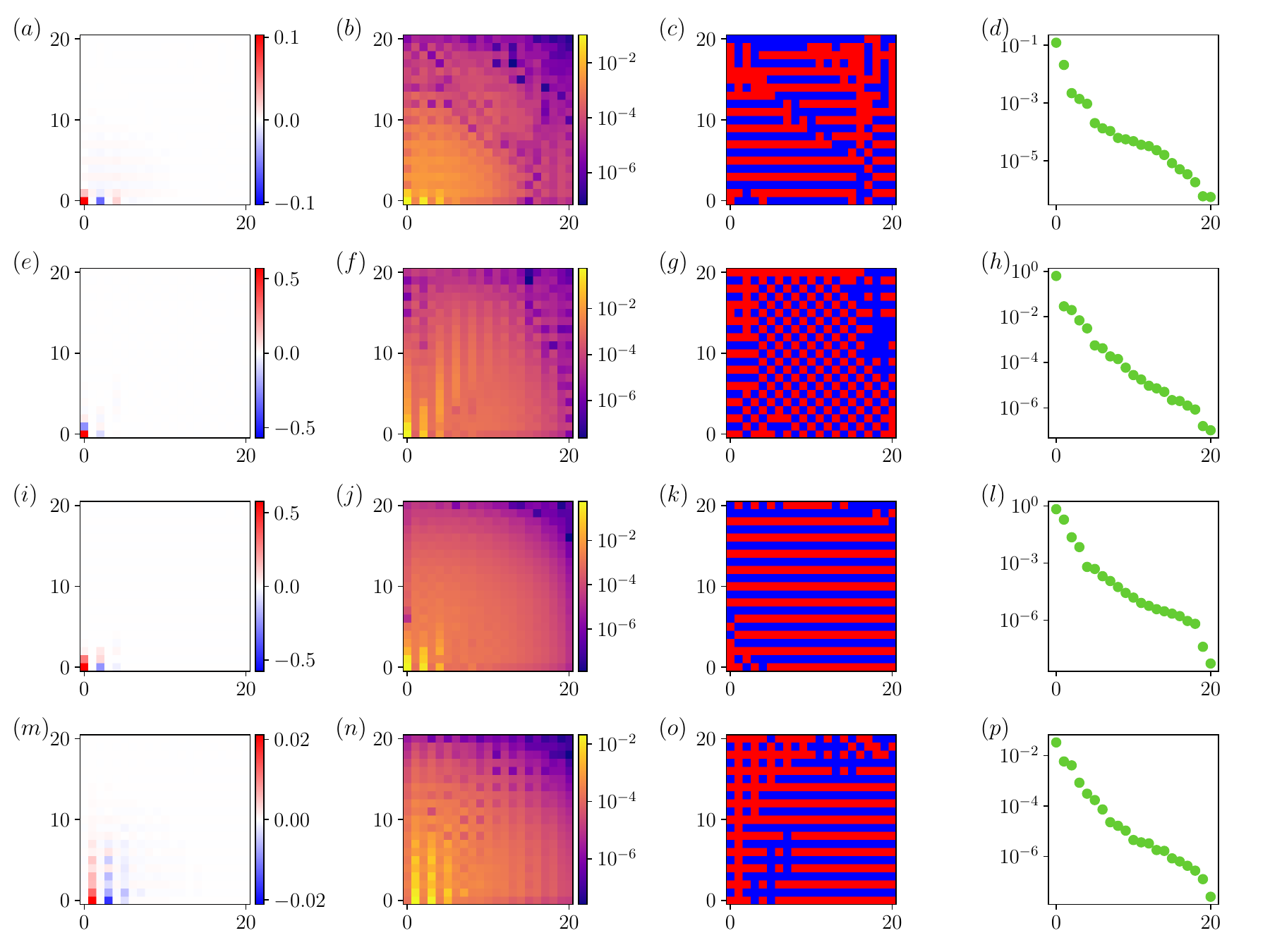}
  \caption{Coefficient matrix for the viscosity inversion solution with the boundary condition shown in the third column of Fig.~\ref{fig:SISyntheticViscosity}. The coefficients were found using a Levenberg Marquardt solver. The top row shows the coefficients for $\mu$, the second for $h$, the third for $u$, and the bottom for $v$. The first column shows the raw coefficient values, the second column the absolute value of the coefficients plotted on a log scale, and the third column the sign of the coefficients. Red indicates positive values while blue represents negative ones. The fourth columns shows the singular values of the coefficient matrices.}
  \label{fig:SIViscCoeffsB}
\end{figure}
\begin{table}
\centering
\caption{Model Hyperparameters and Configuration.}
\label{tab:my_new_table}
\begin{tabular}{@{}lllcc l c cccc@{}}
\toprule
% New row for the spanning header
& & & & & & & \multicolumn{4}{c}{\textbf{Collocation Points}} \\
\cmidrule(lr){8-11} % Adds a partial line under the spanning header
% Main header row with shortened labels
\textbf{Figure} & \textbf{Solver} & \textbf{Loss} & \textbf{Steps} & \textbf{Learning Rate} & \textbf{Basis} & \textbf{Basis Size} & \textbf{PDE} & \textbf{Boundary} & \textbf{Initial} & \textbf{Data} \\
\midrule
2a & Dogleg & Unweighted & --- & --- & Chebyshev & (61,61) & 20,000 & 3,000 & --- & --- \\
2d & LM & Unweighted & --- & --- & Chebyshev & (101,101) & 35,000 & 35,000 & --- & --- \\
3a & Dogleg & Unweighted & --- & --- & CosineLegendre & (201,11) & 30,000 & --- & 1,000 & --- \\
3d & Dogleg & Unweighted & --- & --- & CosineChebyshev & (46,46) & 20,000 & --- & 150 & --- \\
4 & NAdam & Reg. Adap. & 5000 & 0.1 & Legendre & (13,13,13) & 30,000 & 10,000 & 30,000 & --- \\
5 & Adam & Reg. Adap. & 1000 & 0.01 & Chebyshev & (11,11,11) & 15,000 & --- & 4,000 & --- \\
6a & Dogleg & Unweighted & --- & --- & Chebyshev & (26,26) & 20,000 & 1,000 & --- & 20,000\\
6b & Dogleg & Unweighted & --- & --- & Chebyshev & (51,51) & 20,000 & 1,000 & --- & 20,000\\
7 & Adam & Reg. Adap. & 500 & 0.1 & Chebyshev/Fourier & (11,11,11)/(11,11) & 30,000 & 15,000 & 2,000 & 2,000 \\
8 & Adam & Reg. Adap. & 1000 & 0.1 & Chebyshev & (16,16,16) & 30,000 & 30,000 & 10,000 & 10,000\\
S1 & LM & Unweighted & --- & --- & Chebyshev & (21,21) & 2,000 & 600 & --- & 2,000\\
\bottomrule
\end{tabular}
\end{table}
\end{document}